\documentclass{amsart}
\usepackage{amssymb,amsmath,epsfig,yhmath,psfrag}
\usepackage[all]{xy}

\theoremstyle{plain}
\newtheorem{thm}[equation]{Theorem}
\newtheorem{prop}[equation]{Proposition}
\newtheorem{lemma}[equation]{Lemma}
\newtheorem{cor}[equation]{Corollary}
\theoremstyle{definition}
\newtheorem{remark}[equation]{Remark}
\newtheorem{example}[equation]{Example}
\newtheorem{defn}[equation]{Definition}

\numberwithin{equation}{section}

\newcommand{\R}{{\mathbb R}}
\newcommand{\Z}{{\mathbb Z}}
\newcommand{\CC}{{\mathcal C}}

\newcommand{\NN}{{\mathcal N}}
\newcommand{\PP}{{\mathcal P}}
\newcommand{\Pc}{\PP\, \check{}}
\newcommand{\Ct}{C_{\theta}}
\newcommand{\Mt}{M_{\theta}}
\renewcommand{\Re}{\mathop{\rm Re}\nolimits}
\renewcommand{\Im}{\mathop{\rm Im}\nolimits}
\newcommand{\arc}[1]{\stackrel{\frown}{#1}}

\newcommand{\NCA}{NC^{(2)}(A_{n-1})'}

\begin{document}
\title[Polynomials, meanders, and noncrossing partitions]{Polynomials, meanders, and paths in the lattice of noncrossing partitions}
\author{David Savitt}

\address{David Savitt, Department of Mathematics, University of Arizona, 617 N.\
Santa Rita Ave., Tucson, AZ 85721} \email{savitt@math.arizona.edu}

\begin{abstract}  For every polynomial $f$ of degree $n$ with no double roots, there is an
associated family $\mathcal{C}(f)$ of harmonic algebraic curves,
fibred over the circle, with at most $n-1$ singular fibres.  We
study the combinatorial topology of~$\mathcal{C}(f)$ in the
generic case when there are exactly $n-1$ singular fibres.  In
this case, the topology of~$\mathcal{C}(f)$ is determined by the
data of an $n$-tuple of noncrossing matchings on the set
$\{0,1,\ldots,2n-1\}$ with certain extra properties.  We prove
that there are $2(2n)^{n-2}$ such $n$-tuples, and that all of them
arise from the topology of $\mathcal{C}(f)$ for some polynomial
$f$.
\end{abstract}

\maketitle

\section{Introduction}

\begin{figure}
\begin{center}
\resizebox{1.15in}{1.15in}{\includegraphics{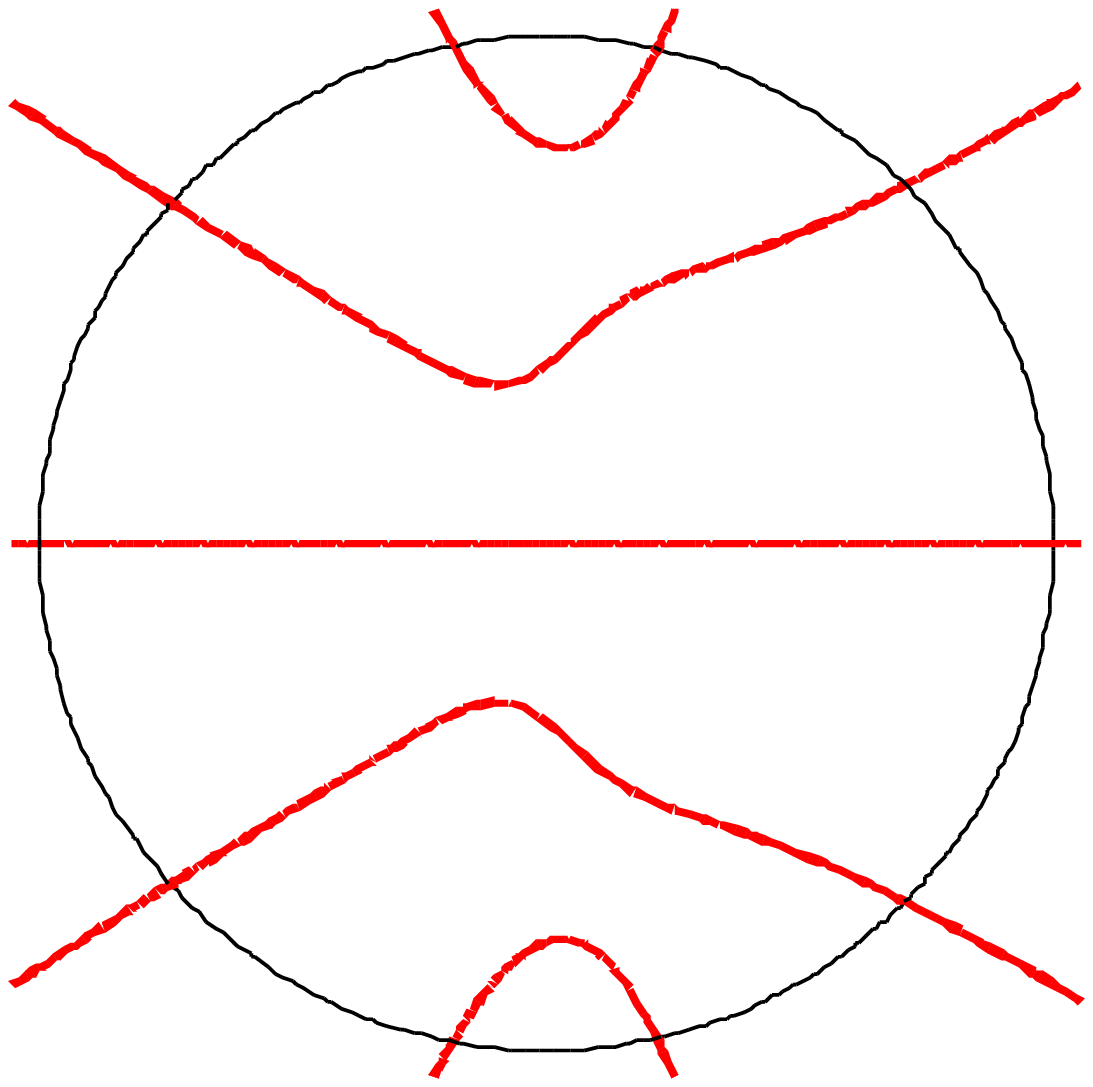}} \hfill
\resizebox{1.15in}{1.15in}{\includegraphics{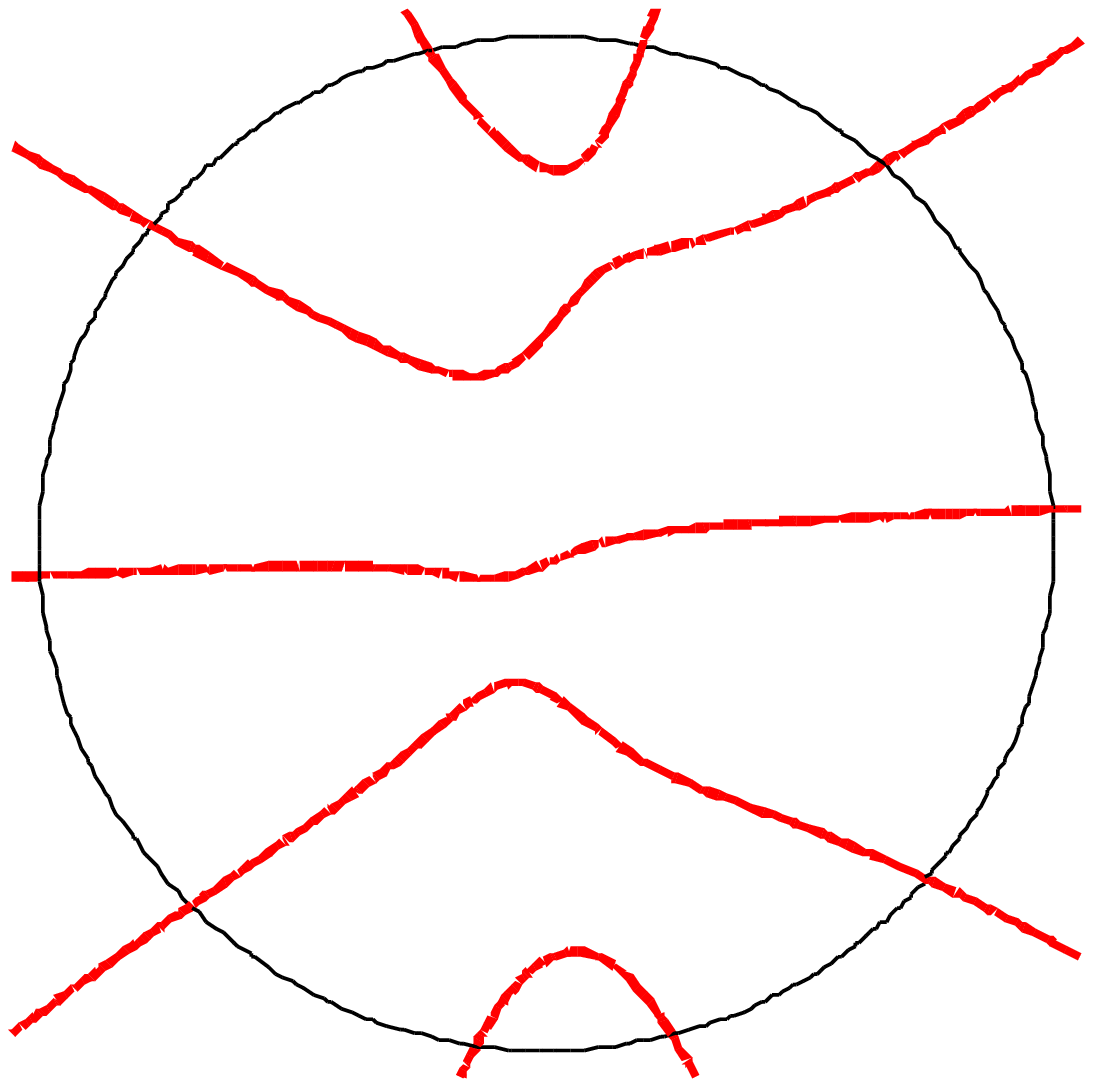}} \hfill
\resizebox{1.15in}{1.15in}{\includegraphics{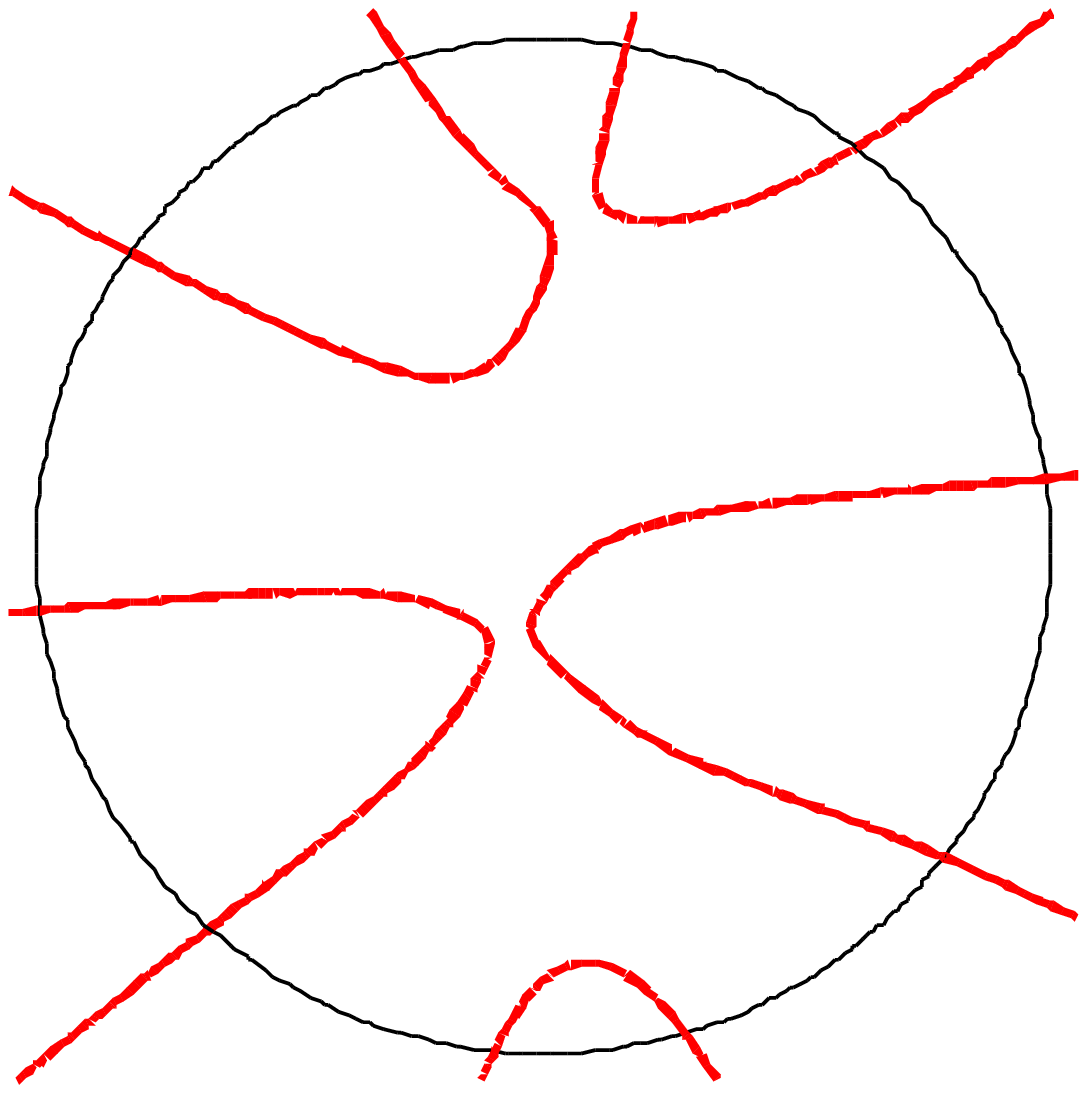}} \hfill
\resizebox{1.15in}{1.15in}{\includegraphics{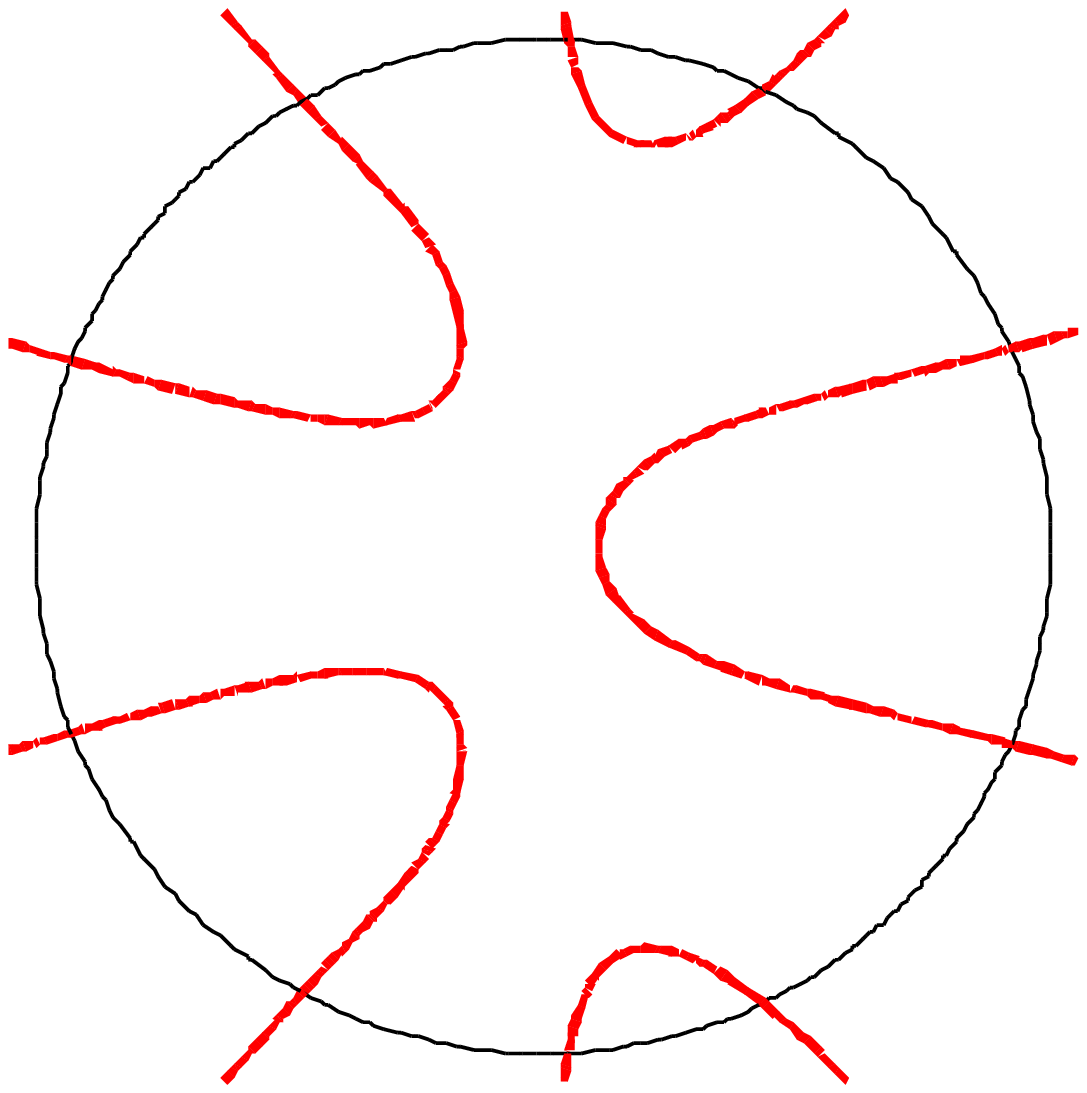}}
\caption{The curves $\Ct(f)$, where $f(z)$ is the quintic $
z^5+6z^3+3z^2+5z-2$, and $\theta=0,
\frac{\pi}{12},\frac{\pi}{6},\frac{\pi}{2}$ respectively.
\label{fig:example-curves}}
\end{center}
\end{figure}

Let $f(z)$ be a monic polynomial of degree $n$ with no double
roots.  For each $\theta \in \R/\pi\Z$, we can associate to $f(z)$
a plane algebraic curve $$\Ct(f) = \{ z \ : \
\Im(e^{-i\theta}f(z)) = 0\}.$$  For instance, $C_{\pi/2}(f) = \{ z
\ : \ \Re(f(z))=0\}$; see Figure \ref{fig:example-curves} for a
pictorial example. Let $\CC(f)$ denote the family of curves
$\Ct(f)$, fibered over the base $\R/\pi\Z$.  Motivated by Gauss's
first proof of the fundamental theorem of algebra, Martin, Singer,
and the author \cite{HarmonicCurves} initiated a study of the
combinatorial topology of the families $\mathcal{C}(f)$, which we
continue in this article.

We begin by recalling the basic properties of the family
$\mathcal{C}(f)$; see \cite{HarmonicCurves} for details. Since
$\Im(e^{-i\theta}f(x+iy))$ is a harmonic function of the variables
$x,y$, the maximum principle dictates that the curves
$C_{\theta}(f)$ do not have any bounded connected components.  In
fact, the curve $C_{\theta}(f)$ has $2n$ asymptotes, at angles
$\frac{\pi k + \theta}{n}$ for integers $0 \le k \le 2n-1$.

The fiber $C_{\theta}(f)$ is singular if and only if there is a root
$r$ of the derivative $f'(z)$ lying on $C_{\theta}(f)$.  By
hypothesis $f(r) \neq 0$ when $r$ is a root of $f'(z)$, so each root
of $f'(z)$ lies on a unique $C_{\theta}(f)$; therefore
$C_{\theta}(f)$ is singular for at most $n-1$ values of $\theta$
(modulo $\pi$).  If $C_{\theta}(f)$ is nonsingular, it has $n$
connected components, each homeomorphic to $\R$, and each containing
two of the $2n$ asymptotes.

Recall that a \textit{matching} on a set of size $2n$ is a
partition of the set into $n$ subsets of size $2$.  Let $\langle
a, b \rangle$ denote the set of integers in the interval $[a,b]$.
If $C_{\theta}(f)$ is nonsingular, then the preceding remarks show
that $C_{\theta}(f)$ induces a matching $M_{\theta}(f)$ on the set
$\langle 0,2n-1 \rangle$, in which $k,k'$ are matched if and only
if the asymptotes at angles $\frac{\pi k + \theta}{n}$ and
$\frac{\pi k' + \theta}{n}$ lie on the same component of
$C_{\theta}(f)$; to make this definition we must choose coset
representatives for $\R/\pi\Z$, and we take $\theta \in [0,\pi)$.
For instance, the matching $M_0(f)$ for the quintic polynomial $f$
of Figure \ref{fig:example-curves} is
$\{\{0,5\},\{1,4\},\{2,3\},\{6,9\},\{7,8\}\}$.

Since the components of $\Ct(f)$ do not cross, $\Mt(f)$ is not
just a matching---it is even a noncrossing matching:

\begin{defn} Let $\PP$ be a partition of a subset of $\R$.
Two blocks of $\PP$ are said to \textit{cross} if there are
integers $i<j<k<\ell$ such that $i,k$ belong to one block and
$j,\ell$ belong to the other block. If no two blocks cross, then
the partition is said to be \textit{noncrossing}.  A
\textit{noncrossing matching (of order $n$)} is a noncrossing
partition (of a set of size $2n$) into blocks of size $2$.
\end{defn}

\begin{figure}[h]
\begin{center}
\resizebox{2.5in}{1in}{\includegraphics{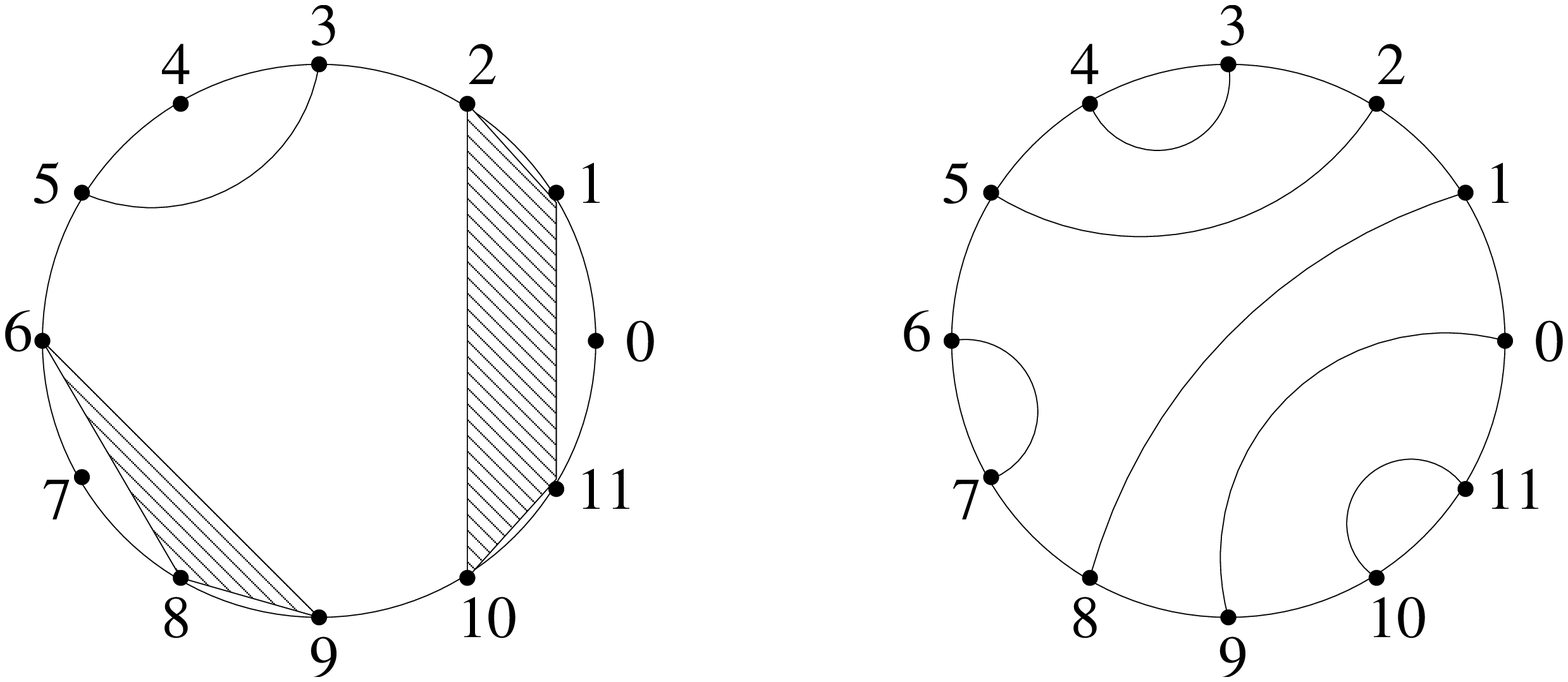}}
\caption{\emph{Left}: The noncrossing partition whose
non-singleton blocks are $\{1,2,10,11\}$, $\{3,5\}$ and
$\{6,8,9\}$.  \emph{Right}: A noncrossing matching.
\label{fig:noncrossing}}
\end{center}
\end{figure}

This definition has a simple geometric interpretation.  Given a
partition of the set $\langle 0,N\rangle$ one can place the integers
from $0$ to $N$ cyclically around a circle, and take the convex
hulls of each of the blocks; the partition is noncrossing if and
only if the convex hulls do not meet.  See Figure
\ref{fig:noncrossing} for examples.

\begin{figure}[h]
\begin{center}
\resizebox{4in}{1in}{\includegraphics{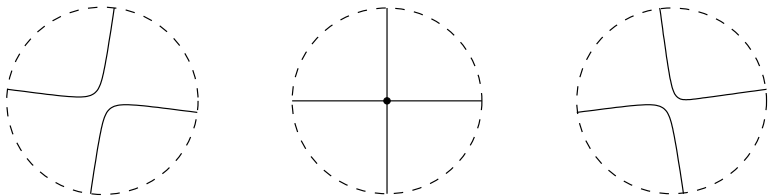}} \caption{The
curves $C_t(f)$ near $z=r$ for $t < \theta$, $t = \theta$, and $t
> \theta$, where $\Ct(f)$ has an ordinary singularity at $r$.} \label{fig:flip}
\end{center}
\end{figure}

Suppose that the fibre $\Ct(f)$ has exactly one singular point
$r$, and suppose furthermore that $f''(r) \neq 0$, so that the
singularity at $r$ is ordinary (nodal).  As $t$
approaches~$\theta$ from below, two components of $C_t(f)$ unite
at the point $r$. As $t$ increases away from $\theta$, the two
components separate in a perpendicular direction, as illustrated
in Figure \ref{fig:flip}.  (We will prove in Proposition
\ref{prop:flip} that this picture is correct.)  Thus, for
$\epsilon$ sufficiently small, the matching
$M_{\theta+\epsilon}(f)$ is a \textit{flip} of the matching
$M_{\theta-\epsilon}(f)$, in the following sense.

\begin{defn} Let $M$ and $M'$ be two matchings on the same underlying set.
We say that $M'$ is a \textit{flip} of $M$ (equivalently, $M$ is a
flip of $M'$) if all but two of the pairs in $M$ are also paired
in $M'$; that is, if $M'$ can be obtained from $M$ by taking two
pairs $\{i,j\},\{k,\ell\}$ in $M$ and replacing them with
$\{i,\ell\},\{j,k\}$ in $M'$.
\end{defn}

Suppose now that $f$ is generic in the sense that the family
$\CC(f)$ has the maximum $n-1$ singular fibres, which are
necessarily ordinary as above.  As $\theta$ varies from $0$ to
$\pi$, each time we cross a singular fibre the matching $\Mt(f)$
changes by a flip.  In this manner we obtain an $n$-tuple $\NN(f)
= (M_1,\ldots,M_n)$ of noncrossing matchings such that $M_{s+1}$
is a flip of $M_s$ for $1 \le s < n$; moreover, we will see that
$M_n$ is the matching $M_1(-1)$, in which $i-1,j-1$ are matched if
and only if $i,j$ are matched in $M_1$ (with subtraction taken
modulo $2n$). The details of this construction can be found in
Section \ref{sec:variation}. Observe that $\NN(f)$ is a
combinatorial invariant which classifies the topology of the
family $\mathcal{C}(f)$.

\begin{defn} A \textit{necklace (of order $n$)} is an $n$-tuple
$\mathcal{N} = (M_1,\ldots,M_n)$ of noncrossing matchings on the
set $\langle 0,2n-1 \rangle$, such that $M_{s+1}$ is a flip of
$M_s$ for $1 \le s < n$, and $M_n = M_1(-1)$.
\end{defn}

\begin{example} Let $M_1 = \{\{0,1\},\{2,3\},\{4,5\}\}$, $M_2 =
\{\{0,3\},\{1,2\},\{4,5\}\}$, $M_3 = \{\{0,5\},\{1,2\},\{3,4\}\}$.
Then $(M_1,M_2,M_3)$ is a necklace of order 3.
\end{example}

In this article we will prove the following two theorems,
answering the questions raised in Section 4 of
\cite{HarmonicCurves}, and classifying the possibilities for the
topology of the family $\CC(f)$ when $f$ is generic as above.

\begin{thm} \label{thm:main-rough} Let $\mathcal{N}$ be a necklace.  Then there exists a
polynomial $f$ such that $\mathcal{N} = \mathcal{N}(f)$.
\end{thm}

\begin{thm} The number of necklaces of order $n$ is $2(2n)^{n-2}$.
\end{thm}

We outline our proof of these theorems.  To begin, recall the main
result of \cite{HarmonicCurves}: if~$M$ is a noncrossing matching
of order $n$, then there exists a polynomial~$f$ such that $M =
\Mt(f)$. This is proved by induction, as follows.  The noncrossing
matching $M$ contains at least one adjacent pair $\{i, i+1\}$,
which we can assume without loss of generality is $\{2n-2,2n-1\}$.
Let $M'$ be the noncrossing matching on $\langle 0, 2n-3 \rangle$
obtained from $M$ by deleting this pair.  By the induction
assumption $M' = \Mt(g)$ for some polynomial $g$, and we show that
$M = \Mt(g(z)(z-z_0))$ for a suitably chosen $z_0$.

Let $\NN = (M_1,\ldots,M_n)$ be the necklace arising from a
polynomial $f$ of degree $n$. We will see in Section
\ref{sec:frompoly} that $\NN$ must satisfy certain conditions
deriving from the fact that $f$ has exactly $n$ roots.  For
instance, if $r < s$, then in the terminology of \cite[Definition
2.4]{HarmonicCurves}, the bimatching corresponding to the pair
$(M_r,M_s)$ must be a basketball.  Any necklace satisfying these
necessary conditions will be called a \textit{strong necklace}. We
will see in Section \ref{sec:frompoly} that the inductive argument
of \cite{HarmonicCurves} can be applied to strong necklaces, so
that every strong necklace arises from a polynomial. Moreover, we
will show that the set of strong necklaces has size $2(2n)^{n-2}$;
although this is a purely combinatorial statement, it is worth
remarking that we will actually invoke the geometry of polynomials
in the proof (more precisely, in the proof of Corollary
\ref{cor:necklaces}). One of the themes of this article is that
the geometry of polynomials can be used to prove nontrivial
statements about noncrossing matchings, and vice-versa.

It remains to prove the combinatorial fact that every necklace is
actually a strong necklace.  We can reinterpret this statement as
follows.

\begin{defn} Let $G_{ncm,n}$ be the graph whose vertices are
the noncrossing matchings of order $n$, and whose edges are the
pairs $(M,M')$ such that $M'$ is a flip of $M$. (Since $n$ will
always be fixed, we will generally omit the subscript $n$ from
$G_{ncm,n}$.)
\end{defn}

Then a necklace is a path of length $n-1$ from $M$ to $M(-1)$ in
$G_{ncm}$, and we wish to show that every such path actually gives
rise to a strong necklace.  (As we will see in Section
\ref{sec:dictionary}, there is an equivalent reinterpretation in
terms of the lattice $NC(n)$ of noncrossing partitions on the set
$\langle 0,n-1 \rangle$.)  The distance from $M$ to $M(-1)$ in
$G_{ncm}$ turns out to be exactly $n-1$: in fact in Section
\ref{sec:meanders} we will give a formula for the distance from
$M$ to $M'$ in $G_{ncm}$ in terms of the number of connected
components of the \textit{system of meanders} associated to the
pair $M,M'$.

This formula in hand, we return in Section \ref{sec:proof} to the
proof of our main result.  Using the machinery we have developed
in the rest of the article, we construct a map from the set of
necklaces to the set of \textit{maximal chains of $2$-divisible
noncrossing partitions of $\langle 0,2n-1 \rangle$} (see
Definition \ref{def:twodivisible}).  This map is seen to be
injective. However, Edelman \cite{Edelman} has shown that the set
of maximal chains of $2$-divisible noncrossing partitions of
$\langle 0,2n-1 \rangle$ has size $2(2n)^{n-2}$.  By our earlier
enumeration of strong necklaces, it follows that this map is a
bijection and that every necklace is a strong necklace, completing
the proof.  In the final section, we use work of Armstrong \cite{Armstrong} to
construct a map from the set of maximal chains of $2$-divisible noncrossing partitions of $\langle 0,2n-1 \rangle$ to the set of necklaces; we are grateful to the anonymous referee for suggesting this argument to us.  As a consequence we obtain some further enumerative properties of necklaces and basketballs.

\section{Variation of $M_{\theta}(f)$} \label{sec:variation}

In this section we describe how $M_{\theta}(f)$ varies with
$\theta$.  The two main technical results (Lemma
\ref{lemma:varies} and Proposition \ref{prop:flip}) are
intuitively transparent; their proofs are simply a matter of
getting the details in order, and can safely be ignored by the
reader.

If $R > 0$, let $D_R$ denote the disk $\{ |z| \le R\}$, and let
$S_R$ denote its boundary circle $\{ |z| = R\}$.

\begin{lemma} \label{lemma:varies}
For $\theta \in [0,\pi)$, suppose that $C_{\theta}(f)$ has a
smooth connected component (necessarily homeomorphic to $\R$)
which contains the asymptotes at angles $\frac{\pi k + \theta}{n}$
and $\frac{\pi k' + \theta}{n}$.
\begin{enumerate}
\item If $\theta \neq 0$, then $C_{t}(f)$ has a smooth connected
component which contains the asymptotes at angles $\frac{\pi k +
t}{n}$ and $\frac{\pi k' + t}{n}$ for all $t$ sufficiently close
to $\theta$.

\item If $\theta = 0$,  then $C_{t}(f)$ has a smooth connected
component which contains the asymptotes at angles $\frac{\pi k +
t}{n}$ and $\frac{\pi k' + t}{n}$ for all $t > 0$ sufficiently
close to $0$, and which contains the asymptotes at angles
$\frac{\pi (k-1) + t}{n}$ and $\frac{\pi (k'-1) + t}{n}$ for all
$t < \pi$ sufficiently close to $\pi$.
\end{enumerate}
\end{lemma}

\begin{proof}
Choose $R \gg 0$ such that
\begin{itemize}
\item $C_{t}(f) \cap S_R$ contains exactly $2n$ points for all
$t$,

\item for all integers $j \in \langle 0,2n-1 \rangle$, the point
$P_j(t)$ of $C_{t}(f) \cap S_R$ which is closest in argument to
$\frac{\pi j + t}{n}$ lies on the asymptote at angle $\frac{\pi j
+ t}{n}$, and

\item $\arg(P_j(t))$ is within $\frac{\pi}{6n}$ of $\frac{\pi j +
t}{n}$.
\end{itemize}

To see that it is possible to choose $R$ satisfying the first
condition, observe by \cite[Remark 1.7]{HarmonicCurves} that it
suffices to show that the set $$\{ r \ : \ S_r \ \text{is tangent
to} \ C_{t}(f) \ \text{for some} \ t\}$$ is bounded, and to take
$R$ larger than this bound. But $C_t(f)$ is tangent to $S_r$ if
and only if the function $h_{t,r}(\lambda) = \text{Im}(e^{-it}
f(re^{i \lambda})) : \R \rightarrow \R$ has a double root.  This
is easily seen to be possible only if $r$ is bounded in terms of
the coefficients of $f$ (and independently of $t$), since
$h_{t,r}(\lambda)$ is the sum of $r^n \sin(n \lambda -t)$ and a
trigonometric polynomial of degree at most $n-1$ whose
coefficients are bounded by $r^{n-1}$ times the coefficients of
$f$.  Now as $R \rightarrow \infty$ the roots of
$h_{t,R}(\lambda)$ approach $\frac{\pi j + t}{n}$ uniformly in
$t$, so the second and third conditions can also be satisfied.

Now let $C$ denote the component of $C_{\theta}(f)$ containing the
asymptotes at angles $\frac{\pi k + \theta}{n}$ and $\frac{\pi k'
+ \theta}{n}$, and suppose that $\delta$ is less than the distance
between $C \cap D_R$ and $(C_{\theta}(f) \setminus C) \cap D_R$,
and also that $\delta < \frac{2R}{3n}$. Set $h_t(z) =
\text{Im}(e^{-it}f(z))$. By item (1) of \cite[Lemma
3.7]{HarmonicCurves} (and in the notation of that lemma),
$C_{t}(f) \cap D_R \subset N_{\delta/2}(C_{\theta}(f) \cap D_R)$
for $t$ sufficiently close to $\theta$; moreover, by item (2) of
\textit{ibid.}, $C_{t}(f) \cap N_{\delta/2}(C) \cap D_R$ is
non-empty.  Assume also that $t$ is sufficiently close to $\theta$
so that for $t \neq \theta$, $C_{t}(f)$ is smooth.  It follows
from the definition of $\delta$ that $C_t(f)$ has a component
$C_t$ such that $C_t \cap D_R$ has endpoints $Q_t,Q'_t$ where
$\arg(Q_t),\arg(Q'_t)$ differ from $\arg P_k(\theta),\arg
P_{k'}(\theta)$ respectively by at most $\frac{\pi}{2R}
\frac{\delta}{2} < \frac{\pi}{6n}$.

Choose $0 < \epsilon < \min(\frac{\pi}{6},\pi-\theta)$ such that
all $t \in (\theta,\theta + \epsilon)$ are sufficiently close to
$\theta$ as above.  By construction, $\arg(Q_t)$ differs from
$\frac{\pi k + t}{n}$ by at most $3 \cdot \frac{\pi}{6n} =
\frac{\pi}{2n}$. It follows that $Q_t = P_k(t)$, and similarly
that $Q'_t = P_{k'}(t)$. Therefore $C_t$ contains the asymptotes
at angles $\frac{\pi k + t}{n}$ and $\frac{\pi k' + t}{n}$, as
desired.

An identical argument works for an interval of the form $(\theta -
\epsilon, \theta)$ with $0 < \epsilon <
\min(\theta,\frac{\pi}{6})$, except that when $\theta = 0$ we must
consider $t$ lying in an interval $(\pi - \epsilon, \pi)$ instead.
In this exceptional case we see that $\arg(Q_t)$ lies within
$\frac{\pi}{2n}$ of $\frac{\pi (k-1) + t}{n} = \frac{\pi k + (t -
\pi)}{n}$, and similarly for $\arg(Q'_t)$.

\end{proof}
\begin{defn}  Let $M$ be a matching on a subset $S$ of the interval
$[0,2n)$.  If $\epsilon$ is any real number, let $S(\epsilon)$
denote the set $\{ i + \epsilon \ : \ i \in S\}$, with addition
taken modulo $2n$.  Define the \textit{rotation of $M$ by
$\epsilon$}, denoted $M(\epsilon)$, to be the matching on
$S(\epsilon)$ such that $i + \epsilon, j+\epsilon$ are matched in
$M(\epsilon)$ if and only if $i,j$ are matched in $M$.
\end{defn}

\begin{prop} \label{prop:constant}
\begin{enumerate}
\item If $U$ is an open interval in $(0,\pi)$ such that
$C_{\theta}(f)$ is nonsingular for all $\theta \in U$, then
$M_{\theta}(f)$ (as a function of $\theta$) is locally constant on
$U$.

\item If $C_0(f)$ is nonsingular, then for all $\theta > 0$
sufficiently close to $0$ the matching $M_{\theta}(f)$ is
$M_0(f)$, while for all $\theta < 0$ sufficiently close to $0$ the
matching $M_{\theta}(f)$ is $M_0(f)(-1)$.
\end{enumerate}
\end{prop}

\begin{proof} To prove (1), apply Lemma
\ref{lemma:varies}(1) to each component of $C_{\theta}(f)$ in
turn, for each $\theta \in U$.  Similarly, to prove (2), apply
Lemma \ref{lemma:varies}(2) to each component of $C_0(f)$ in turn.
\end{proof}

\begin{defn}
\begin{enumerate}
\item We say that $C_{\theta}(f)$ is \textit{completely ordinary}
if $C_{\theta}(f)$ has exactly one singular point, and it is
ordinary; that is, if there is exactly one root $r$ of $f'(z)$
such that $r \in C_{\theta}(f)$, and $f''(r) \neq 0$.

\item We say that a monic polynomial $f(z)$ of degree $n$ is
\textit{completely generic} if $f(z)$ has no double roots, the
family $\mathcal{C}(f)$ has $n-1$ singular fibres (which are
necessarily completely ordinary), and $C_0(f)$ is nonsingular.
\end{enumerate}
\end{defn}

The requirement that $C_0(f)$ be nonsingular stems from our choice
of the interval $[0,\pi)$ as coset representatives for $\R/\pi\Z$
in the definition of $M_{\theta}(f)$.

\begin{prop} \label{prop:flip}  Suppose that $\theta \in (0,\pi)$ and
$C_{\theta}(f)$ is completely ordinary.  Let $M$ be the matching
$M_{t}(f)$ for $t$ approaching $\theta$ from below, and let $M'$
be the matching $M_t(f)$ for $t$ approaching $\theta$ from above.
Then $M'$ is a flip of $M$.
\end{prop}

\begin{proof} Since $C_{\theta}(f)$ is completely ordinary, it has
one singular connected component, and $n-2$ smooth components.
Applying Lemma \ref{lemma:varies}(1) to each of the $n-2$ smooth
components in turn, we see that the matchings $M$ and $M'$ have
(at least) $n-2$ pairs in common.  To conclude that $M'$ is a flip
of $M$, it remains only to dispense with the possibility that
$M'=M$.

Recall our standing hypothesis that $f(z)$ has no double roots.
 Without loss of generality, let us suppose that $C_{\theta}(f)$
is singular at $z = 0$, with $f(0) = e^{i\theta}$, $f'(0) = 0$, and
$f''(0) \neq 0$. By the existence of uniformizers, we can choose a
holomorphic function $g(w)$ on a disk $D'_r$ of radius $r$ centered
at $0$ (in the complex $w$-plane) such that $(f \circ g)(w) = h(w) =
e^{i \theta} (1 + w^2)$.  Let $S'_r$ denote the boundary of $D'_r$.
 Choose $r$ sufficiently small that $g : D'_r \rightarrow g(D'_r)$
and $g \,|_{S'_r}$ are both homeomorphisms.

Observe that $C_{\theta}(h)$ is completely ordinary, consisting of
the real and imaginary $w$-axes (note that $h$ is not monic).  Let
$A_0,\ldots,A_3$ be (short) closed arcs of $S'_r$ such that the
interior of $A_i$ contains $r e^{\pi i/2}$ for $i=0,\ldots,3$. One
checks directly that for $t$ approaching $\theta$ from above,
$C_{t}(h)$ connects a point on $A_0$ to a point on $A_1$ and a point
on $A_2$ to a point on $A_3$, whereas for $t$ approaching $\theta$
from below, $A_0$ and $A_3$ are paired, as are $A_1$ and $A_2$.

Suppose that the singular component of $C_{\theta}(f)$ has a
smooth half-arc (that is, a subset homeomorphic to a real
half-line) connecting a point on $g(A_i)$ to the asymptote at
angle $\frac{\pi k_i + \theta}{n}$.  By an argument identical to
the proof of Lemma \ref{lemma:varies}, the same is true for all
$t$ sufficiently close to $\theta$.  It follows that $M$ contains
the pairs $\{k_0,k_3\},\{k_1,k_2\}$, while $M'$ contains the pairs
$\{k_0,k_1\},\{k_2,k_3\}$.
\end{proof}

Let $f$ be a completely generic polynomial.  Let $\theta_1 <
\cdots < \theta_{n-1}$ be the values of $\theta$ such that
$C_{\theta}(f)$ is singular, and consider the intervals $I_1 =
(0,\theta_1), I_2 =
(\theta_1,\theta_2),\ldots,I_{n}=(\theta_{n-1},\pi)$. From
Proposition \ref{prop:constant}(1) we know that $M_{\theta}(f)$ is
constant on each of the intervals $I_s$; let $M(f)_s$ denote this
matching. By Proposition \ref{prop:flip} we see that $M(f)_{s+1}$
is a flip of $M(f)_s$ for $1 \le s < n$.  By Proposition
\ref{prop:constant}(2) we see that $M(f)_n = M(f)_1(-1)$.

\begin{defn}  If $f$ is a completely generic polynomial, then $$\mathcal{N}(f) =
(M(f)_1,\ldots,M(f)_n)$$ is a necklace, which we call the necklace
arising from $f$.
\end{defn}

Theorem \ref{thm:main-rough} can now be restated more precisely as
follows.

\begin{thm} \label{thm:main-precise} Let $\mathcal{N}$ be a necklace.  Then there exists a
completely generic polynomial $f$ such that $\mathcal{N} =
\mathcal{N}(f)$.
\end{thm}

\begin{remark} We can extend the definition of $\NN(f)$ to
polynomials $f$ which have $n-1$ distinct singular fibres, but for
which $C_0(f)$ is singular, as follows.  We obtain matchings
$M_1,\ldots,M_{n-1}$ just as in the completely generic case, using
the nonsingular intervals
$(0,\theta_1),\ldots,(\theta_{n-2},\pi)$. Now define
$\NN(f)=(M_1,\ldots,M_{n-1},M_1(-1))$.  To check that this is
actually a necklace, note that if we define $g(z)$ via
$f(z)=e^{-i\delta}g(z e^{i \delta/n})$, we will soon see in
\eqref{eq:delta} that $(M_1,\ldots,M_{n-1},M_1(-1))$ is the
necklace of $g(z)$ for $\delta < 0$ sufficiently small .

Let $\overline{G}_{ncm}$ be the quotient of the graph $G_{ncm}$ in
which $M$ and $M(-1)$ are identified for all noncrossing matchings
$M$. Let $\overline{M}$ denote the image of $M$ in the quotient.
Given a polynomial $f$ such that $\CC(f)$ has $n-1$ distinct
singular fibres, it seems to be more intrinsic to associate to $f$
the $(n-1)$-cycle $\overline{\NN}(f)$ containing the vertices
$\overline{M(f)}_1, \ldots, \overline{M(f)}_{n-1}$ in
$\overline{G}_{ncm}$, instead of the necklace $\NN(f)$: this is
free of the choice of coset representatives for $\R/\pi\Z$, and
does not require any ad hoc correction when $C_0(f)$ is singular.
On the other hand, this does lose information, since
$\overline{\NN}(f)=\overline{\NN}(g)$ whenever
$f(z)=e^{-i\delta}g(z e^{i \delta/n})$.  In any event, from a
combinatorial viewpoint necklaces seem to be good objects to
study, as we shall amply see.

\section{Necklaces arising from polynomials} \label{sec:frompoly}

Part of our aim in \cite{HarmonicCurves} was to determine which
ordered pairs $(M,M')$ of noncrossing matchings can arise as a
pair $(M_{\alpha}(f),M_{\beta}(f))$ for $\alpha , \beta \in
[0,\pi)$ with $\alpha < \beta$.  By an elementary parity argument,
each connected component of $C_{\alpha}(f)$ must cross at least
one connected component of $C_{\beta}(f)$, and the point where
they meet is a root of the polynomial $f(z)$ (see the discussion
in \cite[Section 1]{HarmonicCurves}). Since $f(z)$ has exactly $n$
roots, it follows that each connected component of $C_{\alpha}(f)$
crosses exactly one component of $C_{\beta}(f)$, and vice-versa.
It follows that the matching $M_{\alpha}(f) \cup
M_{\beta}(f)(\frac12)$ on the set $\{ 0,\frac12,1,\ldots,2n -
\frac12 \}$ has exactly $n$ crossings.

\begin{defn}  We say that the ordered pair of noncrossing matchings
$(M,M')$ is a \textit{basketball} if the matching $M \cup
M'(\frac12)$ has exactly $n$ crossings.  (This differs slightly
from the terminology of \cite{HarmonicCurves}, where a pair
$(M,M')$ as above would have been called an ordered pair of
matchings corresponding to a basketball.)
\end{defn}

Evidently, for an ordered pair $(M,M')$ of noncrossing matchings
to arise as a pair $(M_{\alpha}(f),M_{\beta}(f))$ as above, it is
necessary that $(M,M')$ be a basketball; \cite[Theorem
3.1]{HarmonicCurves} states that this is sufficient.

Let $\NN = \NN(f) = (M_1,\ldots,M_n)$ be a necklace arising from a
polynomial $f(z)$.  By the above discussion, the pair of
noncrossing matchings $(M_r,M_s)$ must be a basketball for all $1
\le r < s \le n$.  There is another, slightly less obvious,
property of~$\NN$.

\begin{defn} If $(M,M')$ is an ordered pair of matchings on the
set $\langle 0,2n-1 \rangle$ and the blocks $\{i,j\} \in M$ and
$\{i'+\frac12,j' + \frac12\} \in M'(\frac12)$ cross in $M \cup
M'(\frac12)$, we say that $\{i,j\} \in M$ crosses $\{i',j'\} \in
M'$.
\end{defn}

\begin{lemma} Let $\NN = \NN(f)$ be as above.  Suppose $r<s<t$, and
let $\{i_{\star},j_{\star}\}$ be a pair in $M_{\star}$ for each
$\star \in \{r,s,t\}$.  If $\{i_r,j_r\} \in M_r$ crosses
$\{i_s,j_s\} \in M_s$, and $\{i_r,j_r\} \in M_r$ crosses
$\{i_t,j_t\} \in M_t$, then $\{i_s,j_s\} \in M_s$ crosses
$\{i_t,j_t\} \in M_t$.
\end{lemma}

\begin{proof} Let $\lambda_r < \lambda_s < \lambda_t$ be angles in
$[0,\pi)$ such that $M_{\star} = C_{\lambda_{\star}}(f)$ for
$\star \in \{r,s,t\}$.  Let $C_{\star}$ be the connected component
of $C_{\lambda_{\star}}$ corresponding to the matched pair
$\{i_{\star},j_{\star}\}$.  The claim of the lemma is equivalent
to: if $C_r$ meets $C_s$ and $C_t$, then $C_s$ meets $C_t$.  But
the intersection of $C_r$ with both $C_s$ and $C_t$ must be the
unique root of $f(z)$ lying on $C_r$, and the lemma follows.
\end{proof}

\begin{defn} Let $(M_1,\ldots,M_k)$ be a $k$-tuple of noncrossing
matchings on the set $\langle 0,2n-1\rangle$.  We call
$(M_1,\ldots,M_k)$  a \textit{strong pseudonecklace (of length $k$
and order $n$)} (for lack of a better term!) if it satisfies:
\begin{enumerate}
\item The pair of noncrossing matchings $(M_r,M_s)$ is a
basketball for all $1 \le r < s \le k$;

\item If $r<s<t$ and $\{i_r,j_r\} \in M_r$ crosses $\{i_s,j_s\}
\in M_s$, and $\{i_r,j_r\} \in M_r$ crosses $\{i_t,j_t\} \in M_t$,
then $\{i_s,j_s\} \in M_s$ crosses $\{i_t,j_t\} \in M_t$.
\end{enumerate}
\end{defn}

Note that we make no additional hypotheses about how $M_r$ relates
to $M_{r+1}$, or $M_k$ relates to $M_1$.

\begin{example} The reader may verify that if $(M_1,\ldots,M_k)$ is
a strong pseudonecklace, then so is $(M_k(1),M_1,\ldots,M_{k-1})$.
\end{example}

\begin{defn} A strong pseudonecklace $(M_1,\ldots,M_k)$ is a
\textit{strong necklace of length $k$ (and order $n$)} if $M_{r+1}
\neq M_r$ for $1 \le r < k$.  A \textit{strong necklace} (without
reference to its length) will always mean a strong necklace of
length equal to its order $n$.
\end{defn}

It is not immediately clear that a strong necklace is actually a
necklace (i.e., that $M_{r+1}$ is a flip of $M_r$, and $M_n =
M_1(-1)$); this will be a consequence of the fact that every
strong pseudonecklace arises from a polynomial (Theorem
\ref{thm:frompoly}), whose proof begins with the following
observation.

Let $(M_1,\ldots,M_k)$ be a strong pseudonecklace, and let
$\{i_{1,m},j_{1,m}\}$ be an enumeration of the pairs in $M_1$,
with $1 \le m \le n$.  Let $\{i_{r,m},j_{r,m}\}$ denote the
(unique) pair in $M_r$ which crosses $\{i_{1,m},j_{1,m}\}$.  For
each $m$, let $P_m$ be the subset of $\langle 0, 2nk-1\rangle$
consisting of the integers $k i_{r,m} + (r-1)$ and $k j_{r,m} +
(r-1)$ for $1 \le r \le k$.  It is clear that the $P_m$ are
disjoint and exhaust the set $\langle 0,2nk-1\rangle$, i.e.,
together they form a partition $\PP$ of $\langle 0,2nk-1\rangle$.

\begin{lemma} \label{lemma:blocks}  The map $(M_1,\ldots,M_k)
\mapsto \PP$ defines a bijection between the set of strong
pseudonecklaces of length $k$ and order $n$, and the set of
noncrossing partitions of $\langle 0,2nk-1 \rangle$ into $n$
blocks of size $2k$.
\end{lemma}

\begin{proof}  The lemma is clear for $n=1$, so we suppose $n \ge 2$.
First we check that the partition $\PP$ is
noncrossing.  This is clear from the geometric picture, but we
give a formal proof.  Observe first that the pairs $\{k i_{r,m} +
(r-1),kj_{r,m} + (r-1)\}$ and $\{ k i_{s,m'} + (s-1), kj_{s,m'} +
(s-1)\}$ cross if and only if $m=m'$ (and $r \neq s$).  We refer
to the numbers $k i_{r,m} + (r-1)$ and $k j_{r,m} + (r-1)$ as
\textit{counterparts}.

Let $\alpha,\alpha'$ and $\beta,\beta'$ be distinct pairs of
counterparts in the block $P_m$, and suppose without loss of
generality that $\alpha < \beta < \alpha' < \beta'$.  If $m'\neq
m$ and $\gamma,\delta \in P_{m'}$, then we need to prove that
$\gamma,\delta$ are both lie in the same one of the following four
sets: $(\alpha,\beta)$, $(\beta,\alpha')$, $(\alpha',\beta')$, or
$[0,\alpha) \cup (\beta',2nk-1]$.

Suppose, for instance, that $\gamma \in (\alpha,\beta)$.  Let
$\gamma'$ be the counterpart of $\gamma$.  Then $\gamma' \in
(\alpha,\alpha')$ since $\{\gamma,\gamma'\}$ does not cross
$\{\alpha,\alpha'\}$. Similarly $\gamma' \in [0,\beta) \cup
(\beta',2nk-1]$ since $\{\gamma,\gamma'\}$ does not cross
$\{\beta,\beta'\}$.  Hence $\gamma'$ is in $(\alpha,\beta)$ as
well.  An analogous argument works for $\gamma$ in each of the
other three sets. Now we are done if $\delta = \gamma'$. If not,
the same argument shows that $\delta$ and its counterpart
$\delta'$ both lie in the same one of the four sets above. Since
$\{\gamma,\gamma'\}$ crosses $\{\delta,\delta'\}$, all four of
$\gamma,\gamma',\delta,\delta'$ must lie in the same set, and the
partition $\PP$ is indeed noncrossing.

To show that the map $(M_1,\ldots,M_k) \mapsto \PP$ is a
bijection, we construct an inverse.  Let $\PP$ be a noncrossing
partition of $\langle 0,2nk-1\rangle$ into $n$ blocks of size
$2k$.  Let $P$ be a block of $\PP$, containing elements $\alpha_1
< \cdots < \alpha_{2k}$.  Observe that $\alpha_{i+1} \equiv
\alpha_i + 1 \pmod{2k}$.  Indeed, if $\alpha_i + 1 <
\alpha_{i+1}$, then the set $\langle \alpha_i+1 , \alpha_{i+1} - 1
\rangle$ must be a union of blocks of $\PP$ since $\PP$ is
noncrossing, hence it has size divisible by $2k$.  It follows that
each block $P_m$ contains two numbers of the form $k x + (r-1)$
for each $1 \le r \le k$; let these be $k i_{r,m} + (r-1)$ and $k
j_{r,m} + (r-1)$.  For each $1 \le r \le k$ define $M_r$ to be the
collection of pairs $\{ \{i_{r,m},j_{r,m}\} \ : \ 1 \le m \le n
\}$.  It is easy to see that $(M_1,\ldots,M_k)$ is a strong
pseudonecklace, and that this construction is a two-sided inverse
to the map $(M_1,\ldots,M_k) \mapsto \PP$. \end{proof}

\begin{thm} \label{thm:frompoly}  Suppose $0 \le \lambda_1 < \cdots <
\lambda_k < \pi$, and let $(M_1,\ldots,M_k)$ be a strong
pseudonecklace of length $k$ and order $n$.  Then there exists a
monic polynomial $f(z)$ of degree $n$ such that $M_i =
M_{\lambda_i}(f)$ for $1 \le i \le k$.
\end{thm}

\begin{proof}  The statement is trivial for $n=1$.  Assume that the
theorem is true for pseudonecklaces of order $n-1$.

The induction step proceeds in two stages.  First, we show that if
the strong pseudonecklace $(M_1,\ldots,M_k)$ occurs as
$(M_{\lambda_1}(f),\ldots,M_{\lambda_k}(f))$ for all choices of $0
\le \lambda_1 <  \cdots < \lambda_k  < \pi$, then so does the
strong pseudonecklace $(M_k(1),M_1,\ldots,M_{k-1})$.

Let $f(z)$ be a monic polynomial of degree $n$, and define another
monic polynomial $g(z)$ by the formula $f(z) = e^{-i \delta}
g(ze^{i \delta/n})$.  Put $w = z e^{i \delta/n}$.  Then
$\Im(e^{-i\theta} f(z)) = 0$ if and only if
$\Im(e^{-i(\delta+\theta)} g(w)) = 0$, and so the asymptotes at
angles $\frac{\pi k+\theta}{n}$,$\frac{\pi k' +\theta}{n}$ are
matched in $\Ct(f)$ if and only if the asymptotes at angles
$\frac{\pi k+(\theta+\delta)}{n}$,$\frac{\pi k'
+(\theta+\delta)}{n}$ are matched in $C_{\theta+\delta}(g)$.  It
follows that \begin{equation} \label{eq:delta}
M_{\theta+\delta}(g) =
\begin{cases}
\Mt(f) & \text{if} \ 0 \le \theta+\delta < \pi \\
\Mt(f)(1) & \text{if} \ \pi \le \theta+\delta < 2\pi \,.
\end{cases} \end{equation}
Given $0 \le \lambda_1 < \cdots < \lambda_k < \pi$, choose
$\lambda_1 < \delta < \lambda_2$.  By assumption we can find a
monic polynomial $f$ so that
$$(M_{\lambda_2-\delta}(f),\ldots,M_{\lambda_k-\delta}(f),M_{\pi +
\lambda_1 - \delta}(f)) = (M_1,\ldots,M_k)\,.$$ If $g(z)$ is
defined as above, the preceding argument shows that
$$(M_{\lambda_1}(g),\ldots,M_{\lambda_k}(g)) =
(M_k(1),M_1,\ldots,M_{k-1})$$ as desired.  Note that we can use a
similar argument to show that if $(M_1,\ldots,M_k)$ occurs as
$(M_{\lambda_1}(f),\ldots,M_{\lambda_k}(f))$ for all choices of
$\lambda_i$ with $\lambda_1 \neq 0$, then it also occurs for all
choices of $\lambda_i$ with $\lambda_1 = 0$.

For the second stage, note that if $(M_1,\ldots,M_k)$ corresponds
as in Lemma \ref{lemma:blocks} to the partition $\PP$ of $\langle
0,2nk-1\rangle$ into $n$ blocks of size $2k$, then
$(M_k(1),M_1,\ldots,M_{k-1})$ corresponds to the partition
$\PP(1)$.  Observe that any partition $\PP$ of $\langle
0,2nk-1\rangle$ into $n$ blocks of size $2k$ contains at least one
block of the form $\{x,x+1,\ldots,x+2k-1\}$; this is geometrically
clear, and can be proved in an identical manner to
\cite[Proposition 2.11]{HarmonicCurves}.  By the first stage of
the proof, it suffices to prove the statement of the theorem for
strong pseudonecklaces $(M_1,\ldots,M_k)$ corresponding (via Lemma
\ref{lemma:blocks}) to partitions $\PP$ containing the block
$\{0,1,\ldots,k-1,2nk-k,\ldots,2nk-1\}$, i.e., for strong
pseudonecklaces such that $\{0,2n-1\}$ is a pair in $M_i$ for all
$i$.  Moreover, we have seen that we may assume $\lambda_1 > 0$.

But then the result follows directly from \cite[Theorem
3.3]{HarmonicCurves}.  Indeed, let $\check{M}_i$ be the matching
on $\langle 0,2n-3\rangle$ obtained  by omitting the pair
$\{2n-2,2n-1\}$ from $M_i(-1)$.  By the induction assumption
$(\check{M}_1,\ldots,\check{M}_k) =
(M_{\lambda_1}(\check{f}),\ldots,M_{\lambda_k}(\check{f}))$ for
some polynomial $\check{f}$ of degree $n-1$.  Then
$(M_1,\ldots,M_k) = (M_{\lambda_1}(f),\ldots,M_{\lambda_k}(f))$
for $f(z) = \check{f}(z)(z-R)$ for $R \gg 0$.  \end{proof}

\begin{remark} \label{rmk:cg} The proof of Theorem \ref{thm:frompoly} can be
modified slightly to ensure that the polynomial $f$ that one
constructs is completely generic.  Indeed, it suffices to ensure
that $f$ has $n-1$ distinct singular fibres, for one can always
ensure that any particular fibre is nonsingular (in our case,
$C_0(f)$) by replacing $f(z)$ with $g(z)$ given by
$f(z)=e^{-i\delta}g(z e^{i \delta/n})$ for $\delta$ sufficiently
small that $M_{\lambda_i}(g) = M_{\lambda_i-\delta}(f) =
M_{\lambda_i}(f)$, invoking Proposition \ref{prop:constant}(1)
(take $\delta < 0$ if $\lambda_i =0$).  Now by induction
$\check{f}(z)$ can be assumed to be completely generic, and $f$
has distinct singular fibres for $R \gg 0$ by \cite[Lemma
3.5(3),(4)]{HarmonicCurves}.
\end{remark}

\end{remark}

\begin{cor} \label{cor:necklaces} Every strong necklace $\NN = (M_1,\ldots,M_n)$ is a
necklace; that is, $M_{r+1}$ is a flip of $M_r$ for $1 \le r <n$,
and $M_n = M_1(-1)$.  Moreover, $\NN = \NN(f)$ for a completely
generic polynomial $f(z)$.
\end{cor}

\begin{proof}  Pick any $0 < \lambda_1 < \cdots < \lambda_n <
\pi$, and by Theorem \ref{thm:frompoly} choose $f(z)$ monic of
degree~$n$ such that $M_r = M_{\lambda_r}(f)$.  Since $M_r \neq
M_{r+1}$ for $1 \le r <n$, the family $\CC(f)$ must be singular
for some $\theta_r \in (\lambda_r,\lambda_{r+1})$.  Since the
family has at most $n-1$ singularities, the $C_{\theta_r}(f)$ for
$1 \le r \le n-1$ must be all of the singular fibres. Hence $f$ is
completely generic and $(M_1,\ldots,M_n) = \NN(f)$.
\end{proof}

We remark that while the first statement of Corollary
\ref{cor:necklaces} is purely combinatorial, our proof uses the
geometric input of Theorem \ref{thm:frompoly}.  The following two
corollaries have proofs which are similar to that of Corollary
\ref{cor:necklaces}; for Corollary \ref{cor:latter}, invoke Remark
\ref{rmk:cg}.

\begin{cor} There do not exist strong necklaces of length $k$ and
order $n$ with $k > n$.
\end{cor}

\begin{cor} \label{cor:latter} For any strong necklace $(M_1,\ldots,M_k)$ of length $k < n$, there
exists a (strong) necklace $\NN$ such that $(M_1,\ldots,M_k)$ is
obtained by omitting $n-k$ matchings from $\NN$.
\end{cor}

We conclude this section with the following enumerative result.

\begin{prop} \label{prop:enumerate} The number of strong necklaces is $2(2n)^{n-2}$.
\end{prop}

\begin{proof}  Let $\#S$ denote the size of a set $S$, let
$SPN(k,n)$ denote the set of strong pseudonecklaces of length $k$
and order $n$, and let $SN$ denote the set of strong necklaces.
Given a strong pseudonecklace $\NN = (M_1,\ldots,M_k)$ of length
$k$ and order $n$, define $S(\NN) = \{ r < k \ : \ M_{r} = M_{r+1}
\}$. Let $S$ be any subset of $\langle 1,k-1 \rangle$. There is a
bijection between the set $\{ \NN \in SPN(k,n) \ : \ S \subset
S(\NN)\}$ and the set $SPN(k - \#S,n)$, obtained by omitting
$M_{r_i + 1}$ from $\NN$ for each $r_i \in S$.  Note that the
strong necklaces are exactly the strong pseudonecklaces of length
$n$ and order $n$ such that $S(\NN) = \emptyset$.  By the
inclusion-exclusion principle, we have \begin{eqnarray*} \#SN & =
& \sum_{S \subset \langle 1,n-1 \rangle} (-1)^{\#S} \cdot \#\{ \NN
\in SPN(n,n) \ : \ S \subset S(\NN) \} \\
& = &  \sum_{S \subset \langle 1,n-1 \rangle} (-1)^{\#S} \cdot
\#SPN(n-\#S,n) \\
& = & \sum_{j=0}^{n-1} (-1)^j \binom{n-1}{j} \cdot \#SPN(n-j,n)
\,.
\end{eqnarray*}

But by Lemma \ref{lemma:blocks} the strong pseudonecklaces of
length $n-j$ and order $n$ are in bijection with the noncrossing
partitions of $\langle 0,2n(n-j)-1\rangle$ into $n$ blocks of size
$2(n-j)$, and according to \cite[Lemma 4.1]{Edelman} there are
exactly $\frac{1}{n}\binom{2(n-j)n}{n-1}$ of these, so that
$$ \#SN = \frac{1}{n} \sum_{j=0}^{n-1} (-1)^j \binom{n-1}{j} \binom{2(n-j)n}{n-1}
\,.$$

This binomial sum can be evaluated with ease as follows.  Observe
that $\binom{2(n-j)n}{n-1}$ is a polynomial in $j$ of degree $n-1$
with leading coefficient $\frac{1}{(n-1)!} \cdot (-2n)^{n-1}$.  We
have the following general binomial identity \cite[Equation
(5.52)]{GKP}:
$$ \sum_{j=0}^{n-1} (-1)^j \binom{n-1}{j} (a_0 + a_1 j + \cdots
+ a_{n-1} j^{n-1}) = (-1)^{n-1} (n-1)! \cdot a_{n-1} \,.$$
Applying this identity to our sum for $\#SN$ yields
$$ \#SN = \frac{1}{n} (-1)^{n-1} (n-1)! \ \frac{(-2n)^{n-1}}{(n-1)!}
= 2 (2n)^{n-2}$$ as desired.
\end{proof}

\section{The dictionary between $G_{ncm}$ and $NC(n)$}
\label{sec:dictionary}

Before proceeding to our proof that every necklace is a strong
necklace, we pause for a moment to show that the graph $G_{ncm}$
is isomorphic to the graph $NC(n)$ whose definition is as follows.

\begin{defn} The graph $NC(n)$ is the graph
whose vertices are the noncrossing partitions of the set $\langle
0,n-1 \rangle$, and whose edges are pairs $(\PP,\PP')$ of partitions
such that $\PP'$ is obtained from $\PP$ by replacing two blocks of
$\PP$ with their union, or vice-versa.
\end{defn}

In fact $NC(n)$ is a lattice, in which $\PP$ lies below $\PP'$ if
the blocks of $\PP'$ are unions of blocks $\PP$.  It is of course
well-known that the vertices of $G_{ncm}$ and $NC(n)$ are in
bijection, but we wish to record the stronger statement that the
graphs are isomorphic. This observation is not logically necessary
for our proof that every necklace is a strong necklace. However, the
graph $NC(n)$ (or more precisely, the lattice it underlies) has been
well-studied (see e.g. \cite{Kreweras, EdelmanSimion, Simion}) and
so it will be of interest that the properties we eventually
establish for $G_{ncm}$ apply equally well to $NC(n)$.

We recall the standard bijection between the vertices of $G_{ncm}$
and $NC(n)$ (but we omit the proof that the construction is
well-defined and bijective). Let $M$ be a noncrossing matching on
the set $\langle 0,2n-1 \rangle$; it is especially useful here to
think of $M$ as a matching on $2n$ points on a circle, labelled
cyclically from $0$ to $2n-1$. Beginning from $0,1$, etc., relabel
these points $0,0',1,1',\ldots,n-1,(n-1)'$. Now each edge of $M$
joins an unprimed point $a$ to a primed point $b'$, and from the
equivalence relation generated by the relations $a \sim b$ we obtain
a partition $\PP(M)$.

Conversely, suppose we begin with a noncrossing partition $\PP$ on
the set $\langle 0,n-1 \rangle$.  We wish to define a matching
$M(\PP)$.  Suppose that $\{a_{i1},\ldots,a_{i j_i}\}$ is a block
of $\PP$, with $a_{i1} < \cdots < a_{i j_i}$.  We match $a_{i1}'$
with $a_{i2}$, $a_{i2}'$ with $a_{i3}$, and so forth, as well as
$a_{i j_i}'$ with $a_{i1}$.  Having done this for all the blocks,
we obtain a matching on $\{0,0',1,1',\ldots,n-1,(n-1)'\}$. Relabel
these from $0$ to $2n-1$, and we obtain~$M(\PP)$.  The maps $M
\mapsto \PP(M)$ and $\PP \mapsto M(\PP)$ are inverse to one
another.

Next, we recall that the graph $NC(n)$ is self-dual, as follows. If
$\PP$ is a vertex of $NC(n)$, we associate to $\PP$ an element
$\sigma(\PP)$ of the symmetric group $S_n$: write each of the blocks
of $\PP$ as $\{a_{i1},\ldots,a_{i j_i}\}$ with $a_{i1} < \cdots <
a_{i j_i}$, and define $\sigma(\PP)$ to be the product of the cycles
$(a_{i1} \ \cdots \ a_{i j_i})$.  The dual of $\PP$ is the vertex
$\Pc$ satisfying $$ \sigma(\PP) \cdot \sigma(\Pc \,) = (0 \ 1 \
\cdots \ n-1)\,$$ with multiplication taken from right to left.
(Caution: it is not true that $(\Pc \,)\, \check{} = \PP$.)  We will
return to this point of view in Section \ref{sec:enumerative}.  For
details see, e.g., \cite{McCammond} or the memoir \cite{Armstrong}.

The full statement that we wish to establish is as follows.

\begin{prop} The map $M \mapsto \PP(M)$ defines a graph isomorphism
from $G_{ncm}$ to $NC(n)$, under which $\PP(M(-1)) = \PP(M)\,
\check{}$.
\end{prop}

\begin{proof}  If $M$ is a matching
on $\langle 0,2n-1\rangle$, let $\overline{M}$ denote the
corresponding matching on $\{0,0',1,1',\ldots,n-1,(n-1)'\}$. We
prove the duality statement first.   Note that if $\overline{M}$
has an edge from $a'$ to $b$, then $\overline{M(-1)}$ has an edge
from $(b-1)'$ to $a$.  Then $\sigma(\PP(M(-1)))$ sends $b-1$ to
$a$, while $\sigma(\PP(M))$ sends $a$ to $b$, and so
$$\sigma(\PP(M)) \cdot \sigma(\PP(M(-1))) = (0 \ 1 \ \cdots \ n-1);$$
the duality claim follows.

Next, suppose that $\PP'$ is obtained from $\PP$ by replacing two
blocks $\{a_{1},\ldots,a_{j}\}$ and $\{b_1,\ldots,b_k\}$ with their
union, so that $(\PP,\PP')$ is an edge in $NC(n)$.  We wish to show
that $(M(\PP),M(\PP'))$ is an edge in $G_{ncm}$, i.e., that
$M(\PP')$ is a flip of $M(\PP)$ (equivalently, that
$\overline{M(\PP')}$ is a flip of $\overline{M(\PP)}$).  Suppose
without loss of generality that $a_1 < \cdots < a_j$ and $b_1 <
\cdots < b_k$, so that $a_{i}'$ and $a_{i+1}$ are joined by an edge
in $\overline{M(\PP)}$ (with the subscript $i+1$ taken to be $1$ if
$i=j$), and similarly for $b_{i}'$ and $b_{i+1}$.  Without loss of
generality assume $a_1 < b_1$. Since $\PP$ is a noncrossing
partition, the block $\{a_{1},\ldots,a_{j}, b_1,\ldots,b_k\} \in
\PP'$ must have the form $$a_1 < \ldots < a_{\ell} < b_1 < \ldots <
b_k < a_{\ell+1} < \cdots < a_j$$ for some $\ell$ (possibly
$\ell=j$). Then $\overline{M(\PP')}$ is obtained from
$\overline{M(\PP)}$ by flipping the edges $\{a_{\ell}',a_{\ell+1}\}$
and $\{b_k',b_1\}$ to edges $\{a_{\ell}',b_1\}$ and
$\{b_k',a_{\ell+1}\}$ (with $\ell+1$ taken to be $1$ if $\ell=j$).

Conversely, suppose that $M$ is a flip of $M'$, with pairs
$\{a',b\}$ and $\{c',d\}$ in $\overline{M}$ replaced by $\{a',d\}$
and $\{c',b\}$ in $\overline{M'}$.  The only blocks of $\PP(M)$,
$\PP(M')$ affected by this flip are the one or two blocks
containing $a,b,c,d$.  The block of $\PP(M)$ containing $a,b$
corresponds to a sequence of pairs
$\{a',b\},\{b',x_1\},\{x_1',x_2\},\ldots,\{x_j',a\}$ in
$\overline{M}$. If one of these pairs is $\{c',d\} =
\{x_i',x_{i+1}\}$, then the flip breaks this single block into the
two blocks $\{b,x_1,\ldots,x_i=c\}$ and
$\{x_{i+1}=d,\ldots,x_j,a\}$ in $\PP(M')$.  If not, then $c,d$ lie
in another block corresponding to pairs
$\{c',d\},\{d',y_1\},\ldots,\{y_k',c\}$ in $\overline{M}$, and the
flip joins these two blocks into the block
$\{b,x_1,\ldots,x_j,a,d,y_1,\ldots,y_k,c\}$.  In both cases
$(\PP,\PP')$ is an edge in $NC(n)$.  \end{proof}

\begin{example}  We will see shortly that the diameter of the graph $G_{ncm}$ is
$n-1$, and that the distance from $M$ to $M(-1)$ is $n-1$.
Therefore, under the dictionary described in this section, we see
that the diameter of $NC(n)$ is $n-1$, and that the distance from
$\PP$ to $\Pc$ is $n-1$.  Moreover, necklaces correspond to
diameters of the form $\PP \leadsto \Pc$ in $NC(n)$. In
particular, there are exactly $2(2n)^{n-2}$ diameters of this
form.
\end{example}

\section{Meanders and distances in $G_{ncm}$} \label{sec:meanders}

Given two noncrossing matchings $M$ and $M'$ on the set $\langle
0,2n-1\rangle$, we would like to compute the distance between $M$
and $M'$ in the graph $G_{ncm}$.  The key is to consider the
system of meanders associated to the pair $(M,M')$.  Meanders have
been studied previously by Lando and Zvonkin \cite{LandoZvonkin},
Arnol'd, and many others.  Systems of meanders are often defined
geometrically, but in order to remain precise, we give the
following combinatorial definition, which applies equally well to
matchings which are not noncrossing.

\begin{defn} Let $M,M'$ be (not necessarily noncrossing) matchings.
Let $\sim$ be the equivalence relation on the set $\langle
0,2n-1\rangle$ generated by the equivalences $i \sim j$ for every
pair $\{i,j\} \in M$, and $i' \sim j'$ for every $\{i',j'\} \in
M'$.  The \textit{system of meanders} $\Pi_0(M,M')$ is the set of
equivalence classes of $\langle 0,2n-1 \rangle$ under this
relation.  Let $\pi_0(M,M')$ denote the number of equivalence
classes.  If $\pi_0(M,M')=1$ then the pair $(M,M')$ is called a
\textit{meander}.
\end{defn}

\begin{example} We have $\Pi_0(M,M) = M$, and $\pi_0(M,M)=n$.  In fact,
$\pi_0(M,M')=n$ if and only if $M=M'$.
\end{example}

Consider the multigraph $G(M,M')$ with $2n$ vertices labelled from
$0$ to $2n-1$, and with an edge $(i,j)$ for each pair $\{i,j\}$ in
$M$ and $M'$.  The graph is bivalent, and therefore it decomposes
as a collection of loops; each loop alternates between edges
coming from pairs in $M$ and edges coming from pairs in $M'$.  The
blocks of the partition $\Pi_0(M,M')$ are exactly the labels of
each loop, and $\pi_0(M,M')$ is the number of loops, i.e., the
number of connected components of $G(M,M')$. In particular, each
block of $\Pi_0(M,M')$ is equipped with a natural
\textit{unoriented} cyclic ordering: the order obtained by
traversing the corresponding loop.

It is clear from this description that if $M''$ is a flip of $M'$,
then the only loops that are affected by replacing $M'$ with $M''$
in $\Pi_0(M,M')$ are the one or two loops that contain edges
corresponding to the two pairs of $M'$ that are flipped.
Symmetrically, either one or two blocks in $\Pi_0(M,M'')$ are
affected.  Therefore $\pi_0(M,M'') = \pi_0(M,M') + \delta$ with
$\delta \in \{-1,0,1\}$.

\begin{remark} \label{rmk:ncspecial}  If $M,M'$ are noncrossing,
then $\Pi_0(M,M')$ has even more structure: we can orient the
loops of $\Pi_0(M,M')$ in a natural manner, as follows.  Observe
that every edge in $G(M,M')$ has one even endpoint and one odd
endpoint. Since each loop alternates between edges coming from
pairs in $M$ and edges coming from pairs in $M'$, we may orient
every loop so that edges of $M$ flow from even endpoints to odd
endpoints, and edges of $M'$ flow from odd endpoints to even
endpoints.  Let $x \rightarrow y$ denote an oriented edge from $x$
to $y$.

Suppose $M,M'$ are noncrossing, and that the noncrossing matching
$M''$ is obtained from $M'$ by flipping $\{x,y\},\{z,w\}$ to
$\{x,z\}$,$\{y,w\}$. Then $x,w$ have the same parity, and $y,z$
have the same parity; suppose without loss of generality that
$x,w$ are odd. If $x,y,z,w$ are in the same loop $L$ in $G(M,M')$,
then the oriented loop $L$ consists of an edge $x \rightarrow y$,
a path from $y$ to $w$, an edge $w \rightarrow z$, and a path from
$z$ to $x$. Removing the edges $x \rightarrow y$ and $w
\rightarrow z$, and replacing them with edges $x \rightarrow z$
and $w \rightarrow y$ yields two loops. Therefore $\pi_0(M,M'') =
\pi_0(M,M') + 1$.  By a similar argument, if $\{x,y\}$ and
$\{z,w\}$ were in different blocks of $\Pi_0(M,M')$, then the two
blocks become one block in $\Pi_0(M,M'')$, and $\pi_0(M,M'') =
\pi_0(M,M') - 1.$

In summary, if $M,M'$ are noncrossing matchings and the
noncrossing matching $M''$ is a flip of $M'$, then $\pi_0(M,M'') =
\pi_0(M,M') \pm 1$.
\end{remark}

\begin{remark}
Although we shall not make use of this fact, it may be worth
noting that if we define $\{\{0,1\},\ldots,\{2n-2,2n-1\}\}$ to be
an even noncrossing matching, then by Remark \ref{rmk:ncspecial}
each noncrossing matching $M$ obtains a well-defined parity, and
the length of a sequence of flips from $M$ to $M'$ has parity
equal to the difference of the parities of $M$ and $M'$.
\end{remark}

\begin{figure}[t]
\begin{center}
\psfrag{0}{\footnotesize{$0$}} \psfrag{1}{\footnotesize{$1$}}
\psfrag{2}{\footnotesize{$2$}} \psfrag{3}{\footnotesize{$3$}}
\psfrag{4}{\footnotesize{$4$}} \psfrag{5}{\footnotesize{$5$}}
\psfrag{6}{\footnotesize{$6$}} \psfrag{7}{\footnotesize{$7$}}
\psfrag{M}{\footnotesize{$M$}} \psfrag{MM}{\footnotesize{$M'$}}
\resizebox{4in}{1in}{\includegraphics{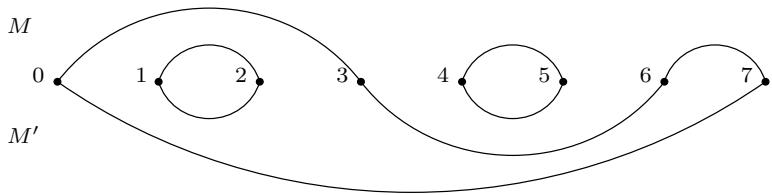}} \caption{A
realization of the system of meanders for the matchings $M$ with
pairs $\{0,3\}$, $\{1,2\}$, $\{4,5\}$, $\{6,7\}$ and $M'$ with
pairs $\{0,7\}$, $\{1,2\}$, $\{3,6\}$, $\{4,5\}$.
 In this example, $\pi_0(M,M')=3$.}\label{fig:meander}
\end{center}
\end{figure}

To study the system of meanders for a pair of noncrossing
matchings, it is helpful to realize the system of meanders
geometrically, as follows.  Consider the points
$\{0,1,\ldots,2n-1\}$ along the real axis in the complex plane.
For each pair $\{x,y\} \in M$, connect the point $x$ on the real
axis to the point $y$ on the real axis via a smooth curve in the
upper half-plane, in such a manner that none of the curves cross;
then do the same in the lower half-plane for the matching $M'$, so
that the whole picture consists of some number of simple closed
curves, as illustrated in Figure \ref{fig:meander}.   The number
of simple closed curves in the picture is $\pi_0(M,M')$, and the
numbers on each curve are the same (and in the same order) as the
numbers on each loop of $G(M,M')$. Such a collection of curves is
called a \textit{realization} of the system of meanders
$\Pi_0(M,M')$.

We will make use of a preferred realization of $\Pi_0(M,M')$,
denoted $C(M,M')$, in which the curve joining each pair is a
semicircle. If $\{x,y\} \in M'$, let $\arc{xy}$ denote the
semicircle joining $x$ and $y$ in the lower half-plane. (We will
use this notation exclusively for pairs in $M'$, and not pairs in
$M$.)  It is not difficult to see that our semicircles do not
cross: if $\{x,y\},\{z,w\}$ are two pairs in $M'$ satisfying $x <
y < z < w$, then no two points on $\arc{xy}$ and $\arc{zw}$ have
the same horizontal coordinate. On the other hand, if $x < z < w <
y$, then $\arc{zw}$ lies above $\arc{xy}$. Up to renaming the
points, these are the only possibilities we need to consider; and
similarly for $M$.  Suppose $x<y$ and $z<w$.  Observe that for any
point $P$ on $\arc{zw}$, the ray extending vertically downwards
from~$P$ intersects $\arc{xy}$ if and only if $x < z < w < y$.

\begin{lemma} \label{lemma:flipexists} If $M,M'$ are noncrossing matchings and
$\pi_0(M,M') < n$, then there exists a noncrossing matching $M''$
such that $M''$ is a flip of $M'$, and $\pi_0(M,M'') = \pi_0(M,M')
+ 1$.
\end{lemma}

\begin{proof}  Since $\pi_0(M,M') < n$, we can find two pairs
$\{x,y\}$ and $\{z,w\}$ in $M'$ which lie in the same block of
$\Pi_0(M,M')$. Without loss of generality we may assume $x < z < w
< y$.  (Up to relabelling, the only other possibility is $x < y <
z < w$; then replace $M,M'$ by $M(2n-w),M'(2n-w)$ to reach the
desired situation.)   We may suppose further that the pairs
$\{x,y\}$ and $\{z,w\}$ have been chosen to minimize the sum
$(z-x) + (y-w)$.   We will prove that the matching $M''$ obtained
by flipping $\{x,y\},\{z,w\}$ to $\{x,z\},\{y,w\}$ is noncrossing.
By Remark \ref{rmk:ncspecial}, this suffices to prove the lemma.

If $\{u,v\}$ crosses $\{x,z\}$ in $M''$, exactly one of the
endpoints of $\arc{uv}$ (let's say $u$) lies in $\langle
x+1,z-1\rangle$; since $\{u,v\}\in M'$ as well, and $M'$ is
noncrossing, $v$ must lie in $\langle w+1, y-1\rangle$.  Therefore,
we wish to show that $M'$ has no pairs with one endpoint in $\langle
x+1,z-1\rangle$ and another in $\langle w+1, y-1\rangle$.  Certainly
no such pair lies in in the block of $\Pi_0(M,M')$ containing
$x,y,z,w$,  by the minimality property of $\{x,y\}$ and $\{z,w\}$.

Let $C$ denote the connected component of $C(M,M')$ containing
$x,y,z,w$.  Pick a point $P$ on $\arc{zw}$, and extend a vertical
ray downwards from $P$ until the ray meets $\arc{xy}$; let $Q$ be
the point of intersection, and let $\overline{PQ}$ denote the line
segment joining $P$ and $Q$.  The interior of $\overline{PQ}$ cannot
meet $C$: by the sentence immediately preceding the statement of the
lemma, any semicircle $\arc{uv}$ with $u<v$ that meets the interior
of $\overline{PQ}$ must have $u \in \langle x+1,z-1\rangle$ and $v
\in \langle w+1, y-1\rangle$, which we have just seen cannot happen
for a semicircle of $M'$ in $C$. Our setup is illustrated in Figure
\ref{fig:flipproof}.

\begin{figure}[bh]
\begin{center}
\psfrag{x}{\footnotesize{$x$}} \psfrag{y}{\footnotesize{$y$}}
\psfrag{z}{\footnotesize{$z$}} \psfrag{w}{\footnotesize{$w$}}
\psfrag{u}{\footnotesize{$u$}} \psfrag{v}{\footnotesize{$v$}}
\psfrag{P}{\footnotesize{$P$}} \psfrag{Q}{\footnotesize{$Q$}}
\psfrag{D}{\footnotesize{$D$}}
\resizebox{4in}{1in}{\includegraphics{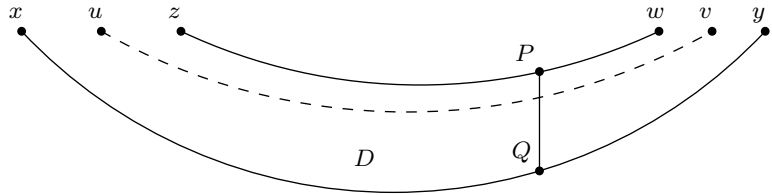}} \caption{The
setup for the proof of Lemma
\ref{lemma:flipexists}}\label{fig:flipproof}
\end{center}
\end{figure}

We imagine that our entire picture is embedded in the sphere (the
one-point compactification of the plane), so that the complement of
$C$ has two connected components, each of which is a topological
open disk. Let $\mathring{D}$ denote the component containing the
interior of $\overline{PQ}$, and let $D$ denote its closure, the
topological closed disk $\mathring{D} \cup C$.  Then $D$ contains
$\overline{PQ}$, and since $\overline{PQ}$ is a path between two
points on the boundary of $D$, it follows that $D \setminus
\overline{PQ}$ has two connected components.

Suppose there exists a pair $\{u,v\}\in M'$ with $u \in \langle
x+1,z-1\rangle$ and $v \in \langle w+1, y-1\rangle$, and let $C'$ be
the connected component of $C(M,M')$ containing $u,v$.  The
semicircle $\arc{uv}$ crosses the interior of $\overline{PQ}$; since
$C'$ does not cross $C$ but does meet the interior of $PQ$, it
follows that $C'$ is contained in $\mathring{D}$.  Since $D
\setminus \overline{PQ}$ has two connected components and $C'$ is a
simple closed curve, it follows that $C'$ meets $\overline{PQ}$ an
even number of times.  (Note that all intersections between
$\overline{PQ}$ and $C(M,M')$ are transverse.)  Therefore $C'$
contains a second semicircle $\arc{rs}$ which meets $\overline{PQ}$,
so that, taking $r<s$, we have $r \in \langle x+1,z-1\rangle$ and $s
\in \langle w+1, y-1\rangle$.  The pairs $\{r,s\}$,$\{u,v\}$
contradict the minimality of the choice of $\{x,y\}$,$\{z,w\}$, so
$u,v$ cannot exist.  This completes the proof. \end{proof}

\begin{defn} Let $G_{mat}$ be the graph whose vertices are the
matchings (not necessarily noncrossing) on $\langle0,2n-1\rangle$,
and whose edges are the pairs $(M,M')$ such that $M'$ is a flip of
$M$.  Let $d_{ncm}$ and $d_{mat}$ denote the distance functions on
$G_{ncm}$ and $G_{mat}$ respectively.
\end{defn}

We are ready to give our formula for distances in the graph
$G_{ncm}$.

\begin{thm} \label{thm:distances} If $M$ and $M'$ are
noncrossing matchings on $\langle 0,2n-1\rangle$, then $$
d_{ncm}(M,M') = d_{mat}(M,M') = n - \pi_0(M,M')\,.$$
\end{thm}

\begin{proof}  Since $G_{ncm}$ is a subgraph of $G_{mat}$, we know
\textit{a priori} that $$d_{mat}(M,M') \le d_{ncm}(M,M')\,.$$  (It
is striking that although $G_{mat}$ is much larger than $G_{ncm}$,
distances between noncrossing matchings will be the same in both
graphs!)

Suppose that $M=M_0,M_1,\ldots,M_k=M'$ is a sequence of (not
necessarily noncrossing) matchings such that $M_{r+1}$ is a flip
of $M_r$ for $r < k$.  Since $\pi_0(M,M)=n$ and $\pi_0(M,M_{r+1})
\ge \pi_0(M,M_r)-1$, we have $\pi_0(M,M') \ge n-k$.  Hence $$n -
\pi_0(M,M') \le d_{mat}(M,M')\,.$$

But by Lemma \ref{lemma:flipexists}, one can produce a path $M' =
M^{0},M^1,\ldots,M^j$ of length $j$ in $G_{ncm}$ such that
\begin{itemize} \item $M^{r+1}$ is a flip of $M^r$ for $r < j$,
\item $\pi_0(M,M^{r+1}) = \pi_0(M,M^r) + 1$ for $r < j$, and \item
$M^j = M$. \end{itemize}
 Hence $n = \pi_0(M,M) = \pi_0(M,M') + j$,
and $j = n - \pi_0(M,M')$.  Therefore
$$ d_{ncm}(M,M') \le n - \pi_0(M,M')\,,$$
and the theorem follows.
\end{proof}

The following corollaries are immediate.

\begin{cor} \label{cor:diameter} The graphs $G_{ncm}$ and $NC(n)$ have diameter $n-1$.
\end{cor}

\begin{cor} \label{cor:meanders} The pair of noncrossing matchings
$(M,M')$ is a meander if and only if
$M,M'$ are a diameter in $G_{ncm}$.
\end{cor}

\begin{remark} Theorem \ref{thm:distances} and Corollaries
\ref{cor:diameter} and \ref{cor:meanders} have been proved
independently by H. Tracy Hall \cite{HTHall}.
\end{remark}

\begin{remark} The proof of Theorem
\ref{thm:distances} shows that $d_{mat}(M,M') = n - \pi_0(M,M')$
for arbitrary matchings $M,M'$ on $\langle 0,2n-1\rangle$.  Note
that the analogue of Lemma~\ref{lemma:flipexists} for arbitrary
matchings is trivial: the difficulty in the proof of
Lemma~\ref{lemma:flipexists} was showing that the matching created
by the desired flip was noncrossing.
\end{remark}

We close this section with the following important observation.

\begin{prop} \label{prop:ismeander} For any noncrossing matching $M$, we have
$\pi_0(M,M(-1))=1$, and $d_{ncm}(M,M(-1))=n-1$.
\end{prop}

\begin{proof} We proceed by induction on $n$; check the claim
directly for $n=1$.

The noncrossing matching $M$ contains at least one pair of the
form $\{i,i+1\}$.  Since $\pi_0(M,M') = \pi_0(M(k),M'(k))$ for any
integer $k$, we can suppose without loss of generality that
$i=2n-2$.  Let $\widehat{M}$ be the noncrossing matching on
$\langle 0,2n-3 \rangle$ obtained by omitting the pair
$\{2n-2,2n-1\}$ from $M$.  Suppose that $\{0,j\} \in M$.

Note that the pairs in $M$ are exactly the same as those in
$\widehat{M}$, with the addition of $\{2n-2,2n-1\}$.  Moreover, the
pairs in $M(-1)$ are exactly the same as those in $\widehat{M}(-1)$,
except for the deletion of $\{j-1,2n-3\}$ and the addition of
$\{j-1,2n-1\}$ and $\{2n-2,2n-3\}$.

By the induction hypothesis, $\pi_0(\widehat{M},\widehat{M}(-1))=1$,
so that any $a,b \in \langle 0,2n-3 \rangle$ are joined by a
sequence of pairs in $\widehat{M}$ and $\widehat{M}(-1)$.  Precisely
the same sequence of pairs in will join them in the system of
meanders of $M$ and $M(-1)$, except that the pair $j-1 \sim 2n-3$
must be replaced by $j-1 \sim 2n-1 \sim 2n-2 \sim 2n-3$ wherever it
occurs. Hence $\langle 0,2n-3\rangle$ is contained in a single block
of $\Pi_0(M,M(-1))$, and since $2n-3 \sim 2n-2 \sim 2n-1$, we
conclude that $\Pi_0(M,M(-1))$ has just one block, as desired.
\end{proof}

\begin{remark} The same statement is false for arbitrary matchings (for
instance, it fails for the matching $M = \{ \{i,n+i\} \ : \ 0 \le
i \le n-1 \}$). \end{remark}

Proposition \ref{prop:ismeander} has the following interesting
geometric consequences.  First, granting Lemma \ref{lemma:varies},
we obtain a new proof of Proposition \ref{prop:flip}, at least in
the case when $f$ is completely generic: from the nonsingular
locus of $\CC(f)$ we obtain $\NN(f)=(M_1,\ldots,M_n=M_1(-1))$ in
which $M_{r+1}$ is either a flip of $M_r$ or equal to $M_r$.  But
if $M_{r+1}=M_r$, we obtain a path from $M_1$ to $M_1(-1)$ in
$G_{ncm}$ of length less than $n-1$, so $M_{r+1}$ must be a flip
of $M_r$.

Second, if $\NN = (M_1,\ldots,M_n)$ is a necklace, then $M_r \neq
M_s$ if $r < s$, or else the sequence of flips $M_1, \ldots, M_r,
M_{s+1}, \ldots, M_n=M_1(-1)$ would be too short a path from $M_1$
to $M_1(-1)$.  (Note that before Proposition \ref{prop:ismeander} we
only knew $M_r \neq M_s$ for $s = r+1$, not arbitrary $s > r$.)
Therefore, if $f(z)$ is a polynomial of degree greater than two, it
is not possible for $\Mt(f)$ to flip from $M$ to $M'$ and then from
$M'$ back to $M$, as~$\theta$ increases.  This does not seem to be
obvious from the geometry alone.

\section{Basketballs and Meanders} \label{sec:proof}

Throughout this section, we will realize basketballs $(M,M')$
geometrically, as follows.  Draw a circle with $4n$ marked points,
labelled $0,\frac12,1,\ldots,2n-\frac12$ counterclockwise around the
circle.   For each pair $\{i,j\}\in M$, draw an arc joining the
points marked $i,j$, and for each pair $\{i,j\} \in M'$, draw an arc
joining the points marked $i+\frac12,j+\frac12$, in such a manner
that each arc of $M$ crosses no arcs of $M$ and exactly one arc of
$M'$, and vice-versa.

To prove that all necklaces are strong necklaces, we must first
begin to understand the systems of meanders associated to
basketballs. We start with some basic properties.

\begin{prop} \label{prop:baskmeanders}
Suppose that the pair of noncrossing matchings $(M,M')$ is a
basketball.  Then:
\begin{enumerate}
\item $\Pi_0(M,M')$ is a noncrossing partition;

\item If the pair $\{s,t\}\in M$ crosses the pair $\{u,v\} \in
M'$, then $s,t,u,v$ lie in the same block of $\Pi_0(M,M')$;

\item If the pair $\{s,t\}\in M$ crosses the pair $\{u,v\} \in M'$
and $s < u + \frac12 < t < v + \frac12$, then the three sets
$\langle s, u\rangle \cup \langle t,v \rangle$, $\langle u+1, t-1
\rangle$, and $\langle v+1,2n-1\rangle \cup \langle 0,s-1\rangle$
are unions of blocks of $\Pi_0(M,M')$.
\end{enumerate}
\end{prop}

\begin{remark} Under the hypotheses of Proposition
\ref{prop:baskmeanders}(3), note that \ref{prop:baskmeanders}(2) and
(3) imply that the block of $\Pi_0(M,M')$ that contains $s,t,u,v$ is
contained in  $\langle s, u\rangle \cup \langle t,v \rangle$.
\end{remark}

\begin{remark} \label{rmk:baskmeanders}
Note that Proposition \ref{prop:baskmeanders}(3) can be restated in
a more uniform manner as follows.  In our geometric realization of
the basketball $(M,M')$, suppose the pair $\{s,t\} \in M$ crosses
the pair $\{u,v\} \in M'$, and the points $s,  u+ \frac12, t,
v+\frac12$ appear in counterclockwise order around the circle.  Then
the following three sets are unions of blocks of $\Pi_0(M,M')$: the
integers which appear counterclockwise between $u + 1$ and $t - 1$;
the integers which appear counterclockwise between $v+1$ and $s-1$;
and the integers which appear counterclockwise between $s$ and $u$
or $t$ and $v$.
\end{remark}

\begin{proof}[Proof of Proposition \ref{prop:baskmeanders}] (2) If one of $s,t$ is equal to one of $u,v$ then
the statement is clear, so we may assume that $s,t,u,v$ are all
distinct.  Replacing $M,M'$ by $M(k),M'(k)$ for a suitable integer
$k$ if necessary, we may suppose without loss of generality that $s
< u + \frac12 < t < v + \frac12$.  Suppose that $t$ is paired with
$x$ in $M'$.  Then $x \not\in \langle s, t \rangle$, or else
$\{s,t\} \in M$ would cross $\{x,t\} \in M'$, and we would have
$v=t$.  Similarly, if $s$ is paired with $y$ in $M'$, we must have
$y \in \langle s+1, t-1 \rangle$.  The loop in $G(M,M')$ containing
$s,t$ contains edges $y - s - t - x$, with $x \not\in \langle s,t
\rangle$ and $y \in\langle s+1,t-1\rangle$.  In order for the loop
to close, it must contain another edge $z - w$ with $z \not\in
\langle s,t \rangle$ and $w \in\langle s+1,t-1\rangle$ (since the
labels $s,t$ cannot occur twice in the loop).  The pair $\{z,w\}$
cannot be a pair in $M$, for it would cross $\{s,t\}$. But then
$\{s,t\} \in M$ crosses $\{z,w\}\in M'$, and since $(M,M')$ is a
basketball, we must have $\{z,w\} = \{u,v\}$.

(1) Let $P,P'$ be distinct blocks of $\Pi_0(M,M')$ which cross, and
suppose $s,t \in P$ and $u,v \in P'$ satisfy $s < u < t < v$. In the
loop of $G(M,M')$ corresponding to the block $P$, there must be an
edge $s' - t'$ such that $s' < u < t' < v$ or $u < t' < v < s'$. In
the former case, in the loop corresponding to $P'$ there must be an
edge $u' - v'$ such that $s' < u' < t' < v'$ or $v' < s' < u' < t'$,
and similarly in the latter case. Since $M,M'$ are noncrossing, one
of $\{s',t'\}$ and $\{u',v'\}$ must be a pair in $M$, and the other
in $M'$. By (2), $s',t',u',v'$ lie in the same block of
$\Pi_0(M,M')$, so $P=P'$.

(3) By Remark \ref{rmk:baskmeanders}, it suffices to prove that
$\langle u+1, t-1 \rangle$ is a union of blocks of $\Pi_0(M,M')$.
Suppose $x \in \langle u+1,t-1\rangle$ and $y \not\in \langle
u+1,t-1 \rangle$ are a pair in either $M$ or~$M'$.  We will derive a
contradiction in each case; the idea is shown in Figure
\ref{fig:proofM}.

\begin{figure}
\begin{center}
\psfrag{x}{\footnotesize{$x$}}
\psfrag{xx}{\footnotesize{$x+\frac12$}}
\psfrag{y}{\footnotesize{$y$?}}
\psfrag{yy}{\footnotesize{$y+\frac12$?}}
\psfrag{s}{\footnotesize{$s$}} \psfrag{t}{\footnotesize{$t$}}
\psfrag{u}{\footnotesize{$u+\frac12$}}
\psfrag{v}{\footnotesize{$v+\frac12$}}
{\includegraphics{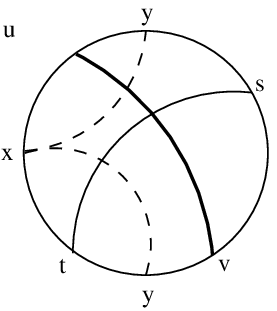}} \qquad \qquad
{\includegraphics{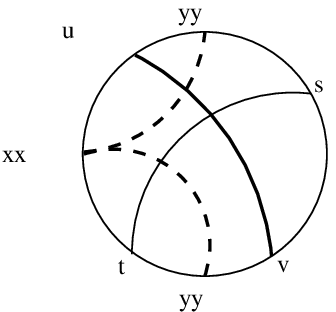}} \caption{The proof of Proposition
\ref{prop:baskmeanders}(3) for $M$ (thin) and $M'$ (\textbf{thick}).
\textit{Left:} the case $\{x,y\} \in M$. \textit{Right:} the case
$\{x,y\} \in M'$. \label{fig:proofM}}
\end{center}
\end{figure}

If $\{x,y\}$ is a pair in $M$, then either $y \in \langle s+1,u
\rangle$ or $y \in \langle t+1, 2n-1\rangle \cup \langle
0,s-1\rangle$ (note that $\{x,y\} \neq \{s,t\}$).  In the former
case, we have $y < u + \frac12 < x < v + \frac12$, contradicting
the hypothesis that only one edge of $M$ crosses $\{u,v\}\in M'$.
In the latter case, we have $y < s < x < t$ or $s < x < t < y$,
contradicting the hypothesis that $M$ is noncrossing.

Similarly, if $\{x,y\}$ is a pair in $M'$, then either $y \in
\langle t,v-1 \rangle$ or $y \in \langle v+1, 2n-1 \rangle \cup
\langle 0, u-1\rangle$ (note that $\{x,y\} \neq \{u,v\}$).  In the
former case, we have $s < x + \frac12 < t < y + \frac12$,
contradicting the hypothesis that $\{s,t\} \in M$ crosses only one
edge in $M'$.  In the latter case, we have $y < u < x < v$ or $u <
x < v < y$, contradicting the hypothesis that $M'$ is noncrossing.
\end{proof}

\begin{remark} \label{rmk:divpartition}  If $M$, $M'$, and $M''$ are
noncrossing matchings and $M''$ is a flip of $M'$, then the
discussion in Remark \ref{rmk:ncspecial} shows that $\Pi_0(M,M'')$
is obtained  by replacing two blocks $\Pi_0(M,M')$ with their
union, or vice-versa.  If   $(M,M')$ and $(M,M'')$ are moreover
both basketballs, then by Proposition \ref{prop:baskmeanders}(1),
$(\Pi_0(M,M'),\Pi_0(M,M''))$ is an edge in $NC(2n)$.
\end{remark}

We can now give the following characterization of basketballs in
terms of systems of meanders.

\begin{prop} \label{prop:bballs} Let $(M,M')$ be a pair of noncrossing matchings.
Then $(M,M')$ is a basketball if and only if \begin{equation}
\label{eq:char} \pi_0(M,M') + \pi_0(M(-1),M') = n+1 \,,
\end{equation}
or equivalently, if and only if $M'$ lies on a path of length
$n-1$ from $M$ to $M(-1)$ in $G_{ncm}$.
\end{prop}

\begin{proof} To see that the two statements are indeed
equivalent, note that \eqref{eq:char} is equivalent to
\begin{equation}
\label{eq:dists} d_{ncm}(M,M') + d_{ncm}(M(-1),M') = n-1
\end{equation} by Theorem \ref{thm:distances}.  Since
$d_{ncm}(M,M(-1))=n-1$, the only way for $M'$ to lie on a path of
length $n-1$ from $M$ to $M(-1)$ in $G_{ncm}$ is if
\eqref{eq:dists} is satisfied.

In the `only if' direction, if $(M,M')$ is a basketball and $M'
\neq M$, then we can apply Corollary \ref{cor:latter} to obtain a
path of length $n-1$ from $M$ to $M(-1)$ in $G_{ncm}$ that
contains $M'$.  If $M'=M$ the `only if' direction is clear.

We will prove the `if' direction by induction on $k =
\pi_0(M,M')$.  If $k=1$, then \eqref{eq:char} implies that
$\pi_0(M(-1),M')=n$, i.e., $M'=M(-1)$.  Since $(M,M(-1))$ is a
basketball, the base case follows.

Now suppose that the statement is known for $k > 0$, and suppose
$M,M'$ satisfy \eqref{eq:char} with $\pi_0(M,M')=k+1$.  Then
$\pi_0(M(-1),M') = n-k < n$, and by Lemma \ref{lemma:flipexists}
there exists a flip $M''$ of $M'$ such that
$\pi_0(M(-1),M'')=n-k+1$.  By the triangle inequality
$$d_{ncm}(M,M'') + d_{ncm}(M(-1),M'') \ge d_{ncm}(M,M(-1))=n-1\,,$$
and so we must have
$$\pi_0(M,M'') + \pi_0(M(-1),M'') \le n+1\,.$$
Since $\pi_0(M,M'') \ge \pi_0(M,M')-1 = k$ by Remark
\ref{rmk:ncspecial}, we conclude $\pi_0(M,M'') = k$.  By the
induction hypothesis $(M,M'')$ is a basketball.  We are therefore
reduced to the following proposition.  (Note that the roles of
$M'$ and $M''$ are reversed in the statement.)
\end{proof}

\begin{prop} \label{prop:basketproof} If $(M,M')$ is a basketball and the noncrossing
matching $M''$ is a flip of $M'$ such that $\pi_0(M,M'') >
\pi_0(M,M')$, then $(M,M'')$ is also a basketball.
\end{prop}

\begin{proof}  Let $\{x,y\},\{z,w\} \in M'$ be pairs
that are flipped to form $\{x,w\},\{y,z\} \in M''$, and without loss
of generality suppose that $x < y < z < w$.  Since $\pi_0(M,M'') >
\pi_0(M,M')$, we know that $x,y,z,w$ are contained in a single block
of $\Pi_0(M,M')$.  What we have to show is that our new pairs $\{
x,w \},\{ y,z \} \in M''$ are each crossed by at most one (and hence
exactly one) pair of $M$.

Suppose that $M$ has a pair $\{s,t\}$ with $t \in \langle
y+1,z\rangle$ and $s \in \langle w+1,2n-1\rangle \cup \langle
0,x\rangle$.  There are (up to easy equivalence) only two
topological possibilities for the pair $\{u,v\}$ in $M'$ that is
crossed by $\{s,t\}\in M$, as shown in Figure \ref{fig:baskchar}.
However, the first picture is impossible since $\{y,z\}$ and
$\{x,w\}$ would both cross $\{u,v\}$ in $M''$; the second picture is
impossible because $z,w$ appear counterclockwise on the circle
between $v+1$ and $s-1$, while $x,y$ appear counterclockwise on the
circle between $s$ and $u$, and so by Proposition
\ref{prop:baskmeanders}(3) and Remark \ref{rmk:baskmeanders} we
could not have $x,y,z,w$ in the same block of $\Pi_0(M,M')$.
Therefore $M$ has no such pair.

\begin{figure}
\begin{center}
\psfrag{x}{\footnotesize{$x+\frac12$}}
\psfrag{y}{\footnotesize{$y+\frac12$}}
\psfrag{z}{\footnotesize{$z+\frac12$}}
\psfrag{w}{\footnotesize{$w+\frac12$}}
\psfrag{s}{\footnotesize{$s$}} \psfrag{t}{\footnotesize{$t$}}
\psfrag{u}{\footnotesize{$u+\frac12$}}
\psfrag{v}{\footnotesize{$v+\frac12$}}
{\includegraphics{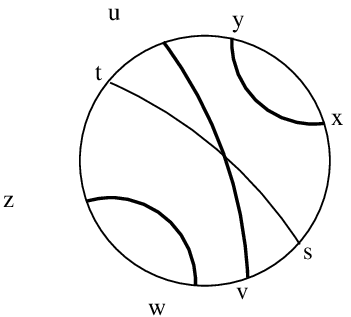}} \qquad \qquad
{\includegraphics{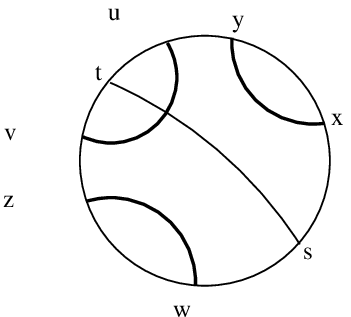}} \caption{The two putative
topological possibilities for $\{u,v\}$ in the proof of Proposition
\ref{prop:basketproof}, with $M$ (thin) and $M'$ (\textbf{thick}).
\label{fig:baskchar}}
\end{center}
\end{figure}

It follows that with the exception of the one pair in $M$ which
crosses $\{x,y\}\in M'$ and the one pair which crosses $\{z,w\} \in
M'$, every other pair of $M$ has both endpoints contained in the
same one of the following four intervals: $\langle z+1,w\rangle$,
$\langle w+1,2n-1\rangle \cup \langle 0,x\rangle$, $\langle
x+1,y\rangle$, and $\langle y+1, z\rangle$.  But no pair with both
endpoints contained in the same one of those four intervals  will
cross either $\{y,z\}\in M''$ or $\{w,x\}\in M''$ (i.e., will cross
$\{ y + \frac12, z+\frac12\}$ or $\{w + \frac12, x+\frac12\}$). Thus
$\{y,z\}\in M''$ and $\{x,w\} \in M''$ are each crossed by at most
the two pairs of $M$ that cross $\{x,y\},\{z,w\} \in M'$; but since
they must be crossed by an odd number of pairs of $M$, they are each
crossed by exactly one pair of $M$.
\end{proof}

\begin{defn} \label{def:twodivisible}  (\cite[Section 4]{Edelman})
A noncrossing partition is called \textit{$k$-divisible} if the
cardinality of every block is divisible by $k$.  On a set of size
$nk$, a \textit{maximal chain of $k$-divisible noncrossing
partitions} is collection of $k$-divisible noncrossing partitions
$(P_1,\ldots,P_{n})$ such that $P_{i}$ has $n+1-i$ blocks and
$P_{i}$ is a refinement of $P_{i+1}$.
\end{defn}

Finally we arrive at our desired conclusion.

\begin{thm} \label{thm:biject}
(1) The map which sends a necklace $(M_1,\ldots,M_n)$ to the
sequence $(\Pi_0(M_1,M_1),\Pi_0(M_1,M_2),\ldots,\Pi_0(M_1,M_n))$ is
a bijection between the set of necklaces of order $n$ and the set of
maximal chains of $2$-divisible noncrossing partitions on $\langle
0,2n-1 \rangle$.

(2) Every necklace is a strong necklace.
\end{thm}

\begin{proof} First we must check that $(\Pi_0(M_1,M_1),\Pi_0(M_1,M_2),\ldots,\Pi_0(M_1,M_n))$
is a maximal chain of $2$-divisible noncrossing partitions on
$\langle 0,2n-1 \rangle$.  Each partition $\Pi_0(M_1,M_i)$ is
noncrossing by Proposition \ref{prop:baskmeanders}(1), and
$2$-divisible since each block is a union of blocks of $M_1$. Since
$\pi_0(M_1,M_1) = n$ and $\pi_0(M_1,M_n) = \pi_0(M_1,M_1(-1))=1$,
and since we know that $\pi_0(M_1,M_{i+1})$ differs by at most $1$
from $\pi_0(M_1,M_i)$, we deduce that $\pi_0(M_1,M_i) = n+1-i$.  By
Remark \ref{rmk:divpartition}, $\Pi_0(M_1,M_i)$ is a refinement of
$\Pi_0(M_1,M_{i+1})$ for each $i < n$. This completes the check.

Next we show that this map is injective. Suppose that a chain of
$2$-divisible noncrossing partitions $P_1,\ldots,P_n$ arises from
a necklace $M_1,\ldots,M_n$; we must show that this necklace is
uniquely determined.  Certainly we must have $M_1 =
\Pi_0(M_1,M_1)=P_1$. Suppose we know that $M_1,\ldots,M_i$ are
determined uniquely, and suppose
$M_1,\ldots,M_i,M_{i+1}',\ldots,M_n'$ is a necklace yielding
$P_1,\ldots,P_n$.

To fix notation, suppose that $P_{i+1} = \Pi_0(M_1,M_{i+1}')$ is
obtained from $\Pi_0(M_1,M_i)$ by replacing the two blocks $B_1,B_2$
with their union.  Returning to the analysis of Remark
\ref{rmk:ncspecial}, observe that in the oriented cyclic ordering on
the block $B_1 \cup B_2$, the elements of $B_1$ are traversed in the
same order as they were in $B_1$ (with some elements of $B_2$
possibly intervening), and similarly for $B_2$. By repeated use of
this observation, it follows that the oriented cyclic ordering on
$B_1 \cup B_2$ in $\Pi_0(M_1,M_{i+1}')$ must be precisely the
oriented cyclic ordering obtained by taking the oriented cyclic
ordering on $\langle 0,2n-1 \rangle$ in $\Pi_0(M_1,M_n')$ and
deleting all elements except those in $B_1 \cup B_2$.  But $M_n' =
M_n = M_1(-1)$, so this oriented cyclic ordering is uniquely
determined by $M_1$!  Given a cyclic ordering on $B_1, B_2$, and
$B_1 \cup B_2$, there is at most one flip of edges in the loops of
$G(M_1,M_i)$ corresponding to $B_1$ and $B_2$ that yields a loop
with the desired ordering on $B_1 \cup B_2$. This must be the flip
which transforms $M_i$ to $M_{i+1}$, and so $M_{i+1}' = M_{i+1}$.
Injectivity follows by induction.

We know that the set of strong necklaces of order $n$ are a subset
of the set of necklaces of order $n$, and we have now shown that
there is an injective map from the set of necklaces of order $n$
to the set of maximal chains of $2$-divisible noncrossing
partitions on $\langle 0,2n-1\rangle$.  From Proposition
\ref{prop:enumerate}, the number of strong necklaces of order $n$
is $2(2n)^{n-2}$.  By \cite[Corollary 4.3]{Edelman}, the number of
maximal chains of $2$-divisible noncrossing partitions on $\langle
0,2n-1\rangle$ is also $2(2n)^{n-2}$.  Both parts of the theorem
follow.
\end{proof}

\section{A map the other way, and enumerative consequences} \label{sec:enumerative}

In this section we construct an injective map from
the set of maximal chains of $2$-divisible noncrossing partitions to
the set of necklaces; combined with the first three paragraphs of
the proof of Theorem \ref{thm:biject}, this furnishes another proof
that these two sets are in bijection.  The
construction relies on recent work of Armstrong \cite{Armstrong}
that was not available when the first version of this article was
written, and we thank the referee for bringing Armstrong's memoir to
our attention, and suggesting the construction. Combined with results from
Section \ref{sec:proof}, the construction yields some further enumerative properties
of basketballs and necklaces.

We begin by setting up some notation.  Let $S_n$ denote the
symmetric group on $\{0,\ldots,n-1\}$, and let $c$ denote the
Coxeter element $(0 \ 1 \ \cdots \ n-1)$.  If $\sigma \in S_n$, its
length $\ell(\sigma)$ is defined to be the least number of
transpositions required to express $\sigma$ as a product of
transpositions.  If $\sigma$ is a product of $r$ disjoint cycles
(including trivial cycles), then $\ell(\sigma)=n-r$.  We write
$\sigma \le \tau$ if $\ell(\sigma) + \ell(\sigma^{-1} \tau) =
\ell(\tau)$, or equivalently if there exists an expression for
$\tau$ as a product of $\ell(\tau)$ transpositions for which the
first $\ell(\sigma)$ transpositions multiply to give $\sigma$. Since
conjugation preserves cycle types, ``first'' can equivalently be
replaced with ``last'' in the preceding sentence.

Let $NC(A_{n-1})$ denote the set of permutations $\sigma$ such that
$\sigma \le c$.  (The notation $NC(A_{n-1})$ is chosen for
consistency with \cite{Armstrong}.)  The following property of
$NC(A_{n-1})$ is standard: if $\PP \in NC(n)$, then the map $\PP
\mapsto \sigma(\PP)$ defined in Section~\ref{sec:dictionary} is a
lattice isomorphism $NC(n) \rightarrow NC(A_{n-1})$, where $\sigma$
lies immediately below $\tau$ in $NC(A_{n-1})$ if there is a
transposition $t$ such that $\tau = \sigma t$ and
$\ell(\tau)=\ell(\sigma)+1$.  If~$M$ is a noncrossing matching on
$\langle 0,2n-1 \rangle$, we let $\sigma(M)$ denote the permutation
$\sigma(\PP(M))$.  Let $\mathcal{S}(n)$ denote the set of $n$-tuples
of permutations $(\sigma(M_1),\ldots,\sigma(M_n))$ associated to
necklaces $(M_1,\ldots,M_n)$.

Let $NC^{(2)}(n)$ denote the set of $2$-divisible noncrossing
partitions of $\langle 0,2n-1\rangle$, and let $\NCA$ denote the
set of all pairs $(\sigma,\tau) \in NC(A_{n-1}) \times
NC(A_{n-1})$ such that $\ell(\sigma \tau) =
\ell(\sigma)+\ell(\tau)$ and $\sigma \tau \le c$.  More
concretely, $(\sigma,\tau) \in \NCA$ if and only if $c$ can be
written as a product of $n-1$ transpositions in which there is a
block of $\ell(\sigma)$ transpositions that multiply to $\sigma$,
to the right of which there is a block of $\ell(\tau)$
transpositions that multiply to $\tau$.  The set $NC^{(2)}(n)$ is
a poset, where $\PP$ lies below $\PP'$ if the blocks of $\PP'$ are
unions of blocks of $\PP$; the set $\NCA$ is a poset, where
$(\sigma,\tau) \le (\sigma',\tau')$ if and only if $\sigma \le
\sigma'$ and $\tau \le \tau'$.  The following is a theorem of
Armstrong, adapted to our notation.

\begin{thm}\cite[Theorem 4.3.8]{Armstrong} \label{thm:arm} There is a poset anti-isomorphism from $NC^{(2)}(n)$
to $\NCA$.
\end{thm}

\begin{remark} To relate Armstrong's statement of this theorem to
the statement above, note that Armstrong's $NC^{(2)}(A_{n-1})$
\cite[Definition 3.3.1(2)]{Armstrong} consists of pairs
$(\pi_1,\pi_2)$ such that $\pi_1 \le \pi_2 \le c$, with the ordering
$(\pi_1,\pi_2) \le (\mu_1,\mu_2)$ if and only if $\mu_1^{-1} \mu_2
\le \pi_1^{-1} \pi_2$ and $\mu_2^{-1}c \le \pi_2^{-1}c$. Then set
$\sigma = \pi_1^{-1} \pi_2$ and $\tau = \pi_2^{-1} c$.
\end{remark}

We now explain how to associate a necklace (more precisely, an
element of $\mathcal{S}(n)$) to a maximal chain of $2$-divisible
noncrossing partitions of $\langle 0,2n-1 \rangle$. It follows
immediately from Theorem \ref{thm:arm} that maximal chains of
$2$-divisible noncrossing partitions are in bijection with
chains  of pairs $(\sigma_1,\tau_1) < \cdots < (\sigma_n,\tau_n)$ in
$\NCA$ such that $\sigma_1 = \tau_1 = 1$, $\sigma_n \tau_n=c$, and
for each $1 \le i < n$ there is a transposition $t_{i}$ such that
either $\sigma_{i+1} = t_i \sigma_i$ and $\tau_{i+1} = \tau_i$, or
else $\tau_{i+1} = t_i \tau_i$ and $\sigma_{i+1} = \sigma_i$.

Let $\ell = \ell(\sigma_n)$ and $\ell' = \ell(\tau_n) = n-1-\ell$.
Observe that there exist $a_1 < \cdots < a_{\ell}$ and $b_1 < \cdots
<  b_{\ell'}$ such that $\sigma_n = t_{a_{\ell}} \cdots t_{a_{1}}$
and $\tau_n = t_{b_{\ell'}} \cdots t_{b_{1}}$.  If $j,j'$ are the
largest integers such that $a_{j},b_{j'} < i$, then $\sigma_i =
t_{a_j} \cdots t_{a_{1}}$ and $\tau_i = t_{b_{j'}} \cdots
t_{b_{1}}$.

To our  chain of pairs $(\sigma_i,\tau_i)$, we associate the
sequence $\rho_i = \sigma_i \tau_n \tau_{i}^{-1}$ for $1 \le i \le
n$. Note that $\rho_1 = \tau_n = \sigma_n^{-1} c$ and
$\rho_n=\sigma_n$. Since the transpositions in the product
\begin{equation} \label{eq:rhoi}
t_{a_j} \cdots t_{a_1} t_{b_{\ell'}} \cdots t_{b_{j'+1}} = \rho_i
\end{equation} are a subsequence of those in the product $\sigma_n \tau_n = c$,
it follows that $\rho_i \le c$ and the product \eqref{eq:rhoi} is a
shortest length product for $\rho_i$. Moreover, $\rho_i$ and
$\rho_{i+1}$ are adjacent in $NC(A_{n-1})$,  with $\ell(\rho_{i+1})
= \ell(\rho_i)+1$ if $i = a_{j+1}$ and $\ell(\rho_{i+1}) =
\ell(\rho_i)-1$ if $i = b_{j'+1}$.  It follows that
$\rho_n,\ldots,\rho_1$ is the sequence of permutations associated to
a necklace from $M$ to $M(-1)$, where $M$ is the matching for which
$\sigma_n = \sigma(M)$ and $\tau_n = \sigma(M(-1))$.

On the other hand, from the sequence $\rho_n,\ldots,\rho_1$ we can
recover the original permutations $(\sigma_i,\tau_i)$ as follows.
Define $\sigma_1 = \tau_1 = 1$.  If $\ell(\rho_{i+1}) =
\ell(\rho_i)+1$ with $\rho_{i+1} = t_i \rho_i$, set $\sigma_{i+1}
= t_i \sigma_i$ and $\tau_{i+1} = \tau_i$. If $\ell(\rho_{i+1}) =
\ell(\rho_i)-1$ with $\rho_i = \rho_{i+1} t_i$, set $\sigma_{i+1}
= \sigma_i$ and $\tau_{i+1} = t_i \tau_i$.  We conclude that the
map sending $(\sigma_1,\tau_1),\ldots,(\sigma_n,\tau_n)$ to
$\rho_n,\ldots,\rho_1$ is one-to-one.  By Theorem
\ref{thm:biject}, it is a bijection onto the set $\mathcal{S}(n)$.

The following is an immediate consequence of the above argument.

\begin{prop} Suppose $M \in G_{ncm}$ and let $\ell = \ell(\sigma(M))$.
Giving a necklace from $M$ to $M(-1)$ is equivalent to giving the
following data: a shortest length factorization of $\sigma(M)$; a
shortest length factorization of $\sigma(M(-1)) = \sigma(M)^{-1} c$;
and a subset of $\langle 1,n-1 \rangle$ of size $\ell$.
\end{prop}

\begin{proof} Given a sequence
$(\sigma_1,\tau_1),\ldots,(\sigma_n,\tau_n)$ corresponding to such a
necklace, so that $\sigma_n = \sigma(M)$, the associated
factorization of $\sigma(M)$ is $t_{a_{\ell}} \cdots t_{a_1}$, the
associated factorization of $\sigma(M)^{-1}c$ is $t_{b_{\ell'}}
\cdots t_{b_1}$, and the associated set of size $\ell$ is
$\{a_1,\ldots,a_{\ell}\}$.  In the reverse direction, the factorizations of $\sigma(M)$ and $\sigma(M)^{-1}c$ determine the sets $\{\sigma_1,\ldots,\sigma_n\}$ and $\{\tau_1,\ldots,\tau_n\}$ in the obvious  manner, and the subset of size $\ell$ in $\langle 1,n-1\rangle$ is the set of integers $i$ such that $\sigma_{i+1}\neq \sigma_i$; its complement is the set of integers $i$ such that $\tau_{i+1} \neq \tau_i$.
\end{proof}

\begin{defn} If $\sigma$ is a permutation of length $\ell$ and cycle type $(m_1,\ldots,m_{n-\ell})$ including trivial cycles, set
$$ \mathcal{B}(\sigma) = \binom{\ell}{m_1-1,\ldots,m_{n-\ell}-1} m_1^{m_1-2} \cdots m_{n-\ell}^{m_{n-\ell}-2} $$
and
$$ \mathcal{C}(\sigma) = C_{m_1} \cdots C_{m_{n-\ell}}$$
where $C_m$ denotes the $m$th Catalan number.
\end{defn}
Then $\mathcal{B}(\sigma)$ is the number of factorizations of $\sigma$ into $\ell$ transpositions: indeed, it is well-known that the number of shortest factorizations of an $m$-cycle is $m^{m-2}$ (see, e.g., \cite{Kreweras}); each minimal factorization of $\sigma$ consists of a minimal factorization of each cycle, and the binomial coefficient gives the number of ways of interleaving the factorizations of the cycles.  Similarly, $\mathcal{C}(\sigma)$ is the number of permutations $u$ such that $u \le \sigma$.

\begin{cor}  Let $\sigma = \sigma(M)$.  The number of necklaces from $M$ to $M(-1)$ is $$\binom{n-1}{\ell(\sigma)} \mathcal{B}(\sigma) \mathcal{B}(\sigma^{-1}c)\,.$$
\end{cor}

\begin{prop} Let $\sigma = \sigma(M)$.  The number of noncrossing matchings $M'$ such that $(M,M')$ is a basketball is
$$ \mathcal{C}(\sigma) \mathcal{C}(\sigma^{-1}c)\,.$$
\end{prop}

\begin{proof}  By Proposition \ref{prop:bballs}, the noncrossing matchings $M'$ such that $(M,M')$ is a basketball are precisely the noncrossing matchings which lie on necklaces from $M$ to $M(-1)$; by the preceding discussion, these are the noncrossing matchings $M'$ such that $\sigma(M') = uv$ with $u \le \sigma$ and $v \le \sigma^{-1}c$.  To complete the proof, we have to check that the product $uv$ uniquely determines $u$ and $v$.

Let $\wedge$ denote meet (greatest lower bound) in the lattice $NC(A_{n-1})$.  We will show that $u = \sigma \wedge uv$; certainly $u \le \sigma \wedge uv$.  Set $K_{\mu}^{\nu}(w) = \mu w^{-1} \nu$.  If $\mu \le \pi \le \nu$, then \cite[Theorem 2.6.14]{Armstrong} gives the formula
$$ K_{\mu}^{\nu}(\pi) = \nu \wedge K_{\mu}^{c}(\pi) \,.$$
Putting $\mu=u$ and $\nu = \pi = \sigma$ gives $u = \sigma \wedge u(\sigma^{-1}c)$.  Then
$$ u \le \sigma \wedge uv \le \sigma \wedge u(\sigma^{-1}c) = u\,,$$
completing the proof.
\end{proof}

\vskip 0.3cm

 \noindent \textbf{Acknowledgments.}  It is a pleasure
to thank Jeremy Martin, Vic Reiner, and Sergei Tabachnikov for
valuable conversations and correspondence, the
Max-Planck-Institut f\"ur Mathematik for its hospitality, and the
anonymous referee for his or her comments. Jeremy Martin graciously
permitted the re-use of Figures \ref{fig:example-curves} and
\ref{fig:noncrossing}, which he drew originally for
\cite{HarmonicCurves}.

\bibliographystyle{math}
\bibliography{savitt}

\end{document}